\documentclass[preprint,review,sort&compress,3p]{elsarticle}
\usepackage[utf8]{inputenc} 
\usepackage{mathrsfs}
\usepackage{amsfonts}
\usepackage{cancel}
\usepackage{graphicx}
\usepackage{amssymb}
\usepackage{hyperref}
\usepackage{amsthm}
\usepackage{amsmath}
\usepackage{epstopdf}
\usepackage{verbatim}
\usepackage{booktabs}%
\usepackage{flafter}
\usepackage{algorithm}
\usepackage{algorithmicx}
\usepackage{algpseudocode}
\usepackage{makecell}
\usepackage{caption,subcaption}
\usepackage{multirow}
\usepackage{appendix}
\usepackage{float}
\usepackage{array, caption, threeparttable}
\usepackage[font=small,labelfont=bf,labelsep=none]{caption}
\numberwithin{equation}{section}

\journal{ }
\begin{document}
	\begin{frontmatter}
		\title{OBK-RCM: Accelerated Orthogonal Block Kaczmarz Algorithm via RCM Reordering and Dynamic Grouping for Sparse Linear Systems}
		\cortext[cor1]{Corresponding author}
		\author[a]{Yu-Fang Liang}
		\author[a]{Hou-Biao Li\corref{cor1}}
		\ead{lihoubiao0189@163.com}
		
		\address[a]{School of Mathematical Sciences, University of Electronic Science and Technology of China, Chengdu, 611731, P. R. China}
		
		\begin{abstract}
  Existing block Kaczmarz methods face challenges in balancing computational efficiency and convergence for large sparse linear systems with scattered nonzero patterns, due to costly partitioning strategies and non-orthogonal projections. In this paper, we propose the orthogonal block Kaczmarz (OBK-RCM) algorithm with the Reverse Cuthill-McKee (RCM), which integrates the RCM reordering with a novel orthogonal block partitioning strategy. RCM transforms sparse matrices into banded structures to enhance inter-block orthogonality, while dynamic grouping of mutually orthogonal blocks based on angle cosine thresholds reduces iterative complexity. In addition, two extended versions (SOBK-RCM and UOBK-RCM) are proposed to deal with non-square systems by constructing extended matrices without sacrificing sparsity. This work offers a practical framework for efficient sparse linear algebra solvers. Experiments on 33 real-world and synthetic matrices show that OBK-RCM achieves 10-50 times faster CPU time (up to several hundred) and 50–90\% fewer iterations than state-of-the-art methods (RBK, RBK(k), GREBK(k), aRBK), especially for scattered sparse structures in most cases. Theoretical analysis confirms linear convergence, driven by hyperplane orthogonality. 
	\end{abstract}
		
	\end{frontmatter}
	
	\section{Introduction}
This paper focuses on solving large-scale sparse linear systems of the form
	\begin{equation}\label{eq:1.1}
		Ax=f,
	\end{equation}
where $A\in {{\mathbb{R}}^{m\times n}}$, $f\in {{\mathbb{R}}^{m}}$, and $x\in {\mathbb{R}^{n}}$. Nowadays, one of the classic and effective iterative projection methods for solving such systems is the Kaczmarz method \cite{1,2}.
	
	The classic Kaczmarz method scans each row of the matrix in a cyclic manner and projects its current iterative solution $x_{j}$ onto a hyperplane $H_{i_k}= \{x\in {\mathbb{R}^{m}}\mid  A_{(i_k)}x = f_{i_k}\}$ defined by the system rows \cite{1,2,20}. Specifically, the initial value is set to $x_{0}$, assuming that the $i_k$th row has been selected at the $j$th iteration, then the $(j+1)$th estimate vector $x_{j+1}$ is obtained by:
	\begin{equation}\label{eq:1.2}
		x_{j+1} = x_{j} + \frac{f_{i_k}-A_{(i_k)}x_j}{\left\|A_{(i_k)}\right\|_{2}^{2}}A_{(i_k)}^T.
	\end{equation}

   The classic Kaczmarz method suffers from slow convergence due to its dependence on row selection order and lacks rigorous theoretical guarantees. To address these limitations, Strohmer and Vershynin \cite{3,4} introduced the randomized Kaczmarz (RK) algorithm, achieving an exponential convergence rate in expectation. However, RK’s random row selection criterion remains suboptimal \cite{4}. Subsequent studies proposed greedy variants: Ansorge developed the greedy residual Kaczmarz (GEK) \cite{5}, Nutini et al. formulated the greedy distance Kaczmarz (GDK) \cite{6}, and Bai and Wu designed the greedy random Kaczmarz (GRK) with relaxation \cite{7,8}. Further extensions include Liu et al.’s GRK adaptation for ridge regression problem \cite{9} and Du et al.’s distance-based GDRK algorithm \cite{10}, which outperforms classical RK and GRK in efficiency.
	
	The block Kaczmarz (BK) method \cite{11,12} is a natural extension of the classical Kaczmarz method. Unlike the row Kaczmarz method, the block Kaczmarz algorithm involves selecting multiple row indices in each iteration. Given an initial solution $x_0$, in the $k$-th iteration, first select a subset of row indices $\tau_j\in\left[m\right]$, and then project the iterative solution $x_{j-1}$ onto the solution space $A_{\tau_j}x = f_{\tau_j}$. Using the pseudo-inverse of $A_{\tau_j}$, the Kaczmarz iteration formula for the block is:
	\begin{equation}\label{eq:1.3}
		x_{j+1} = x_{j} + A_{\tau_j}^{\dagger}\left(f_{\tau_j}-A_{\tau_j}x_{j}\right),
	\end{equation} 
	where $A_{\tau_j}$ and $f_{\tau_j}$ respectively represent the submatrix and right-hand vector corresponding to $A$ and $f$. 

	Building on these advancements, Needell et al. introduced the randomized block Kaczmarz (RBK) algorithm \cite{13,14}. Subsequent work integrated greedy strategies into RBK, yielding the greedy random block Kaczmarz (GRBK) method \cite{17}. Zhang et al. further enhanced residual-based selection with the greedy random Motzkin-Kaczmarz (GRMK) algorithm \cite{15}, while Jiang and Li incorporated K-means clustering to develop the row-clustered RBK(k) \cite{16} and global randomized block variants \cite{17}. Zheng et al. refined this framework by analyzing residual standardization, proposing the greedy residual block Kaczmarz (GREBK(k)) \cite{18}. Most recently, Zhang et al. (2024) combined fast projection and weighted averaging in the aRBK algorithm \cite{19}, achieving a balance between speed and robustness.
	
	While block Kaczmarz methods (e.g., RBK \cite{13}, GREBK(k) \cite{18}, aRBK \cite{19}) improve convergence over row-wise iterations, two critical limitations hinder their efficiency:
	\begin{itemize}
	  \item High partitioning cost: Existing strategies (e.g., K-means clustering [16]) require 
	$O(mk)$ computations for block division, becoming prohibitive for large-scale sparse matrices.
	  \item Slow convergence due to non-orthogonal projections: Random or greedy block selection often results in adjacent blocks with small angles between hyperplanes (Figure 1), leading to redundant iterations.
	\end{itemize}
	
	To address these gaps, we propose the orthogonal block Kaczmarz (OBK-RCM) algorithm with RCM, featuring three core innovations:
	\begin{itemize}
	  \item RCM reordering: The Reverse Cuthill-McKee reordering concentrates nonzeros near the diagonal, forming an approximate banded matrix with enhanced inter-block orthogonality (Section 3.1).
	  \item Dynamic orthogonal partitioning: Blocks with near-zero angle cosines are grouped into orthogonal classes (Oclass), enabling rapid convergence via sequential orthogonal projections (Section 3.3).
	  \item Non-square system compatibility: OBK-RCM extends to handle non-square systems without sacrificing sparsity, through extended matrix construction (SOBK-RCM for $m>n$, UOBK-RCM for $m<n$) (see, Eq. 4.1–4.3), avoiding information loss (Section 4.0).
	\end{itemize}
	
	This work bridges the gap between low-cost partitioning and accelerated convergence, offering a unified framework for sparse linear systems.

\section{Introduction to existing block methods}
	The emergence of the block Kaczmarz method has greatly improved the running time and convergence rate of the Kaczmarz method in handling high-dimensional linear systems. 	
	\begin{algorithm}[htbp]
		\caption{ RBK(k) Algorithm \cite{15}}
		\label{2}
		\begin{algorithmic}[1]  	
			\Require $A$, $f$, $x_{0}$, $k$, $l$, $\theta$
			\Ensure An estimation $x_{l}$ of the unique solution $x_{*}$ to $Ax=f$.
			\State Apply the K-means method to $A$ and $f$ to obtain a partition $\{\tau_1,...,\tau_m\}$ of the row $\{1,...m\}$, and obtain $k$ groups $A_{\tau_i}$ and $f_{\tau_i}$, $i=1,2,...,k$. $\overline{A}=[\overline{A}_{\tau_1},...,\overline{A}_{\tau_k}]$, $\overline{f}=[\overline{f}_{\tau_1},...,\overline{f}_{\tau_k}]$, where $\overline{A}_{\tau_i}$ and $\overline{f}_{\tau_i}$ respectively represent the center of the class.
			\For {$j=0,1,2, \ldots$} until termination criterion is satisfied
			\State Compute
			$\epsilon_j={\frac{\theta}{\left\|\overline{f}-\overline{A}x_j\right\|_{2}^{2}}\mathop{max}\limits_{1 \leq \tau_j \leq k} \left\{\frac{\lvert\overline{f}_{\tau_j}-\overline{A}_{\tau_j}x_j\rvert^{2}}{\left\|\overline{A}_{\tau_j}\right\|_{2}^{2}}\right\}+\frac{1-\theta}{\left\|\overline{A}\right\|_{F}^{2}}},~~~~\theta\in (0,1)$
			\State Define the index set of positive integers $U_{j}=\left\{\tau_j \mid{\lvert\overline{f}_{\tau_j}-\overline{A}_{\tau_j}x_j\rvert^{2}}\geq{\epsilon_j\left\|\overline{f}-\overline{A}x_j\right\|_{2}^{2}\left\|\overline{A}_{\tau_j}\right\|_{2}^{2}}\right\}$
			\State Calculate the $\tau$-th element of vector $\widetilde{r}_j$,
			${\overset{\sim}{r}}_{j}^{(\tau)}=\left\{\begin{aligned}
				\overline{f}_{\tau}-\overline{A}_{\tau}x_j ~~~~if~\tau \in U_j \\
				0~~~~~~~ otherwise \end{aligned}\right.$
			\State Select $\tau_j$ from $U_j$ by probability $P_r(\tau=\tau_j)=\frac{\lvert{\overset{\sim}{r}}_{j}^{(\tau)}\rvert^{2}}{\left\|\widetilde{r}_j\right\|_{2}^{2}}$
			\State $x_{j+1}=x_{j}+A_{\tau_j}^{\dagger}\left(f_{\tau_j}-A_{\tau_j}x_{j}\right)$
			\EndFor
		\end{algorithmic}
	\end{algorithm}

	Firstly, in 2022, Jiang et al. proposed a randomized block Kaczmarz (RBK(k)) method by combining the K-mean clustering algorithm with greedy randomization technique \cite{15}, see Algorithm 1 for details.

	Secondly, in 2024, based on the relationship between common residuals and standardized residuals, Zheng et al. \cite{17} proposed the block Kaczmarz algorithm (GREBK(k)) for clustering and chunking of standardized residuals combined with the K-means algorithm, see Algorithm 2 for details.
	\begin{algorithm}[htbp]
		\caption{GREBK(k) Algorithm \cite{17}}
		\label{3}
		\begin{algorithmic}[1]  	
			\Require $A$, $f$, $x_{0}$, $k$, $l$, $\theta$
			\Ensure An estimation $x_{l}$ of the unique solution $x_{*}$ to $Ax = f$.
			\State Apply the K-means method to the normalized residual vector $d$ to obtain a partition $\{\tau_1,...,\tau_m\}$ of the row $\{1,...m\}$, and obtain $k$ groups $A_{\tau_i}$ and $f_{\tau_i}$, $i=1,2,...,k$. $\overline{A}=[\overline{A}_{\tau_1},...,\overline{A}_{\tau_k}]$, $\overline{f}=[\overline{f}_{\tau_1},...,\overline{f}_{\tau_k}]$, where $\overline{A}_{\tau_i}$ and $\overline{f}_{\tau_i}$ respectively represent the center of the class.
			\For {$j=0,1,2, \ldots$} until termination criterion is satisfied
			\State Compute $\epsilon_j={\frac{\theta}{\left\|\overline{f}-\overline{A}x_j\right\|_{2}^{2}}\mathop{max}\limits_{1 \leq \tau_j \leq k} \left\{\frac{\lvert\overline{f}_{\tau_j}-\overline{A}_{\tau_j}x_j\rvert^{2}}{\left\|\overline{A}_{\tau_j}\right\|_{2}^{2}}\right\}+\frac{1-\theta}{\left\|\overline{A}\right\|_{F}^{2}}},~~~~\theta\in (0,1)$
			\State Define the index set of positive integers $U_{j}=\left\{\tau_j \mid{\lvert\overline{f}_{\tau_j}-\overline{A}_{\tau_j}x_j\rvert^{2}}\geq{\epsilon_j\left\|\overline{f}-\overline{A}x_j\right\|_{2}^{2}\left\|\overline{A}_{\tau_j}\right\|_{2}^{2}}\right\}$
			\State Calculate the $\tau$-th element of vector $\widetilde{r}_j$,
			${\overset{\sim}{r}}_{j}^{(\tau)}=\left\{\begin{aligned}
				\overline{f}_{\tau}-\overline{A}_{\tau}x_j ~~~~if~\tau \in U_j \\
				0~~~~~~~ otherwise \end{aligned}\right.$
			\State Select $\tau_j$ from $U_j$ by probability
			$P_r(\tau=\tau_j)=\frac{\lvert{\overset{\sim}{r}}_{j}^{(\tau)}\rvert^{2}}{\left\|\widetilde{r}_j\right\|_{2}^{2}}$
			\State $x_{j+1}=x_{j}+A_{\tau_j}^{\dagger}\left(f_{\tau_j}-A_{\tau_j}x_{j}\right)$
			\EndFor
		\end{algorithmic}
	\end{algorithm}
	
	In the same year, Zhang et al. \cite{19} proposed a variant of the random Kaczmarz algorithm, the aRBK algorithm, by combining block projection and weight averaging techniques. Their combination can balance convergence speed, convergence range, and computational complexity. See Algorithm 3 for details.
	\begin{algorithm}[htb]
		\caption{aRBK Algorithm \cite{19}}
		\label{4}
		\begin{algorithmic}[1]  	
			\Require $A$, $f$, $x_{0}$, $\left\{\omega_j\right\}_{j\in \tau}$, $\alpha$, $K$, $k=0$
			\Ensure An estimation of the unique solution $x_{*}$ to $Ax = f$.
			\State The application definition divides the system into $k$ groups $A_j$ and $f_j$, $j\in \tau$.
			\While {$k<K$} 
			\State $k=k+1$
			\For {$j\in \tau$}, perform independent sampling $S_j$ for each subsystem, then update
			\State $x_k^j = x_{k-1} + \alpha\left(A_{S_j}^j\right) ^{\dagger}\left(f_{S_j}^j-A_{S_j}^jx_{k-1}\right)$
			\EndFor
			\State Compute $x_k = \sum_{j=1}^{\tau}\omega_jx_k^j$
			\EndWhile
		\end{algorithmic}
	\end{algorithm}

	These algorithms (detailed in Algorithms 1-3) demonstrate how hybrid strategies - clustering of residuals, standardized error metrics, and projection balancing - overcome traditional Kaczmarz limitations in sparse system solving.
	
\section{The orthogonal block Kaczmarz (OBK-RCM) algorithm with RCM}\label{sec3}
\subsection{Matrix reordering to enhance orthogonality}\label{sec3.1}
\begin{figure}[htbp]
		\centering
		\includegraphics[scale=0.45]{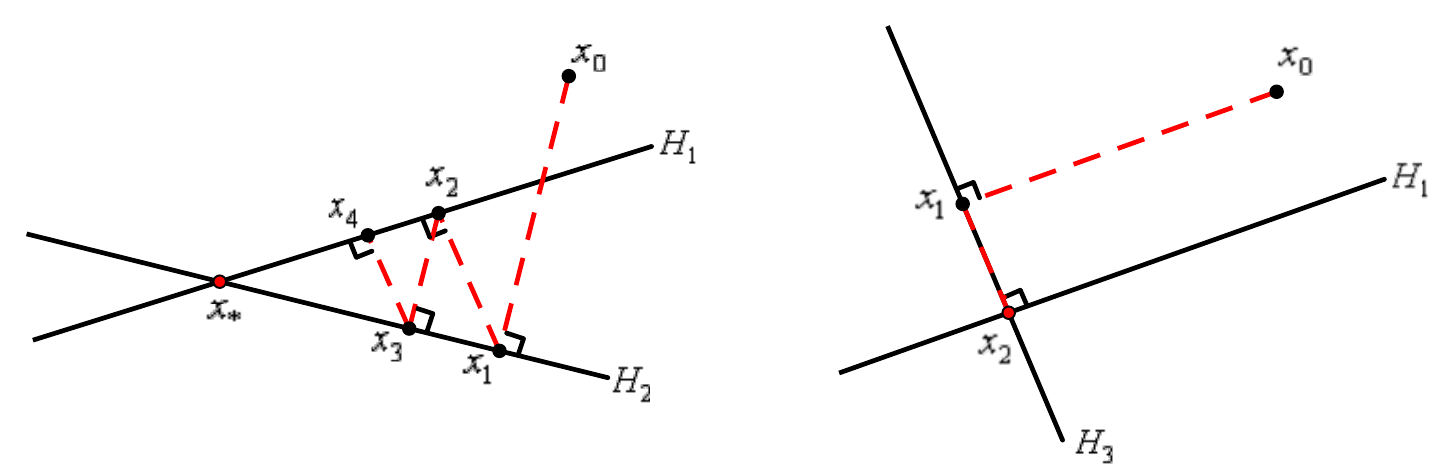}
		\caption{\;Comparison of Kaczmarz Iterations on Small and Vertical Angles on a 2D System.}
		\label{fig1}
	\end{figure}

	The Kaczmarz algorithm exhibits good convergence properties for linear systems, particularly with sparse matrices. However, its convergence rate slows significantly for large-scale sparse matrices or those with scattered structural distributions, even when enhanced by block clustering and greedy algorithms. This limitation arises because the algorithm's convergence speed depends critically on the angle between consecutive projection hyperplanes during iterations. Specifically, smaller angles between hyperplanes necessitate more iterations and slower convergence, whereas near-orthogonal projections (angles close to 90°) yield faster convergence with fewer steps. To illustrate, consider a two-dimensional linear system (Figure \ref{fig1}). Here, $H_1$, $H_2$, $H_3$ denote the system's hyperplanes, $x_0$ the initial point, and $x_1$, $x_2$, $x_3$, $x_4$ the iterative solutions. The red dots where the hyperplanes intersect represent the exact solution $x_*$. The left panel shows slow convergence when $H_1$ and $H_2$ form a small angle, while the right panel demonstrates rapid convergence when $H_1$ and $H_3$ are orthogonal—requiring only two iterations to reach the exact solution $x_*$, compared to four for non-orthogonal cases. This motivates the exploration of methods to extract or construct mutually orthogonal rows in the system matrix, which could drastically reduce iteration time (see \cite{22}).

	It is well known that the row blocks of a banded matrix have good orthogonality between them, and the smaller the bandwidth, the stronger the orthogonality. Based on the above analysis, we propose to process the original matrix into an approximate banded matrix by reordering. Sparse matrix reordering algorithms, such as Reverse Cuthill-McKee (RCM), Column Counting (Colperm), Minimum Degree (Amd), Nested Dissection (Dissect)\cite{23,24,25,26}, are matrix reduction bandwidth efficient processing techniques. Although there are several sorting methods available, experiments have shown that the RCM algorithm tends to provide better computational results in terms of bandwidth reduction, improved inter-block orthogonality, and algorithmic time-consumption. Therefore, for the block Kaczmarz algorithm, the RCM algorithm is preferred. Section 3.2 describes the RCM algorithm and its advantages in detail.

\subsection{Reverse Cuthill-McKee Algorithm}\label{sec3.2}
	The inverse Cuthill-McKee algorithm \cite{23} is a reordering method for bandwidth minimization of sparse matrices. The algorithm drastically reduces the bandwidth of the matrix and enhances the orthogonality of the blocks by rearranging the order of the columns of the matrix rows and clustering the non-zero elements of the original matrix around the diagonal. 

	The RCM, Colperm, Amd, and Dissect algorithms have varying results in reducing bandwidth and enhancing orthogonality. We use the matlab built-in function $sprandn$ to randomly generate sparse matrices of different dimensions. Table \ref{tab1} shows the original bandwidths ($bw$) of these sparse matrices, as well as the bandwidths after reordering by the RCM, Colperm, Amd, and Dissect algorithms. It can be seen that the RCM algorithm can reduce the matrix bandwidth more effectively than the other three sorting algorithms.
	\begin{table}[!hpt]
		\centering
		\setlength {\tabcolsep}{5mm}
		\setlength{\abovecaptionskip}{10pt}%
		\setlength{\belowcaptionskip}{10pt}%
		\caption{\;\;Comparison of processing matrix bandwidth of different sorting algorithms.}
		\begin{tabular}{lcccccc}
			\toprule
			$Name$ && $bw$ & $Colperm$ & $Amd$ & $Dissect$ & $RCM$  \\
			\hline
			$1000\times1000$ && 967 & 984 & 970 & 994 & 28 \\
			$5000\times5000$ && 4959 & 4904 & 4720 & 4837 & 32  \\
			$10000\times10000$ && 9779 & 9813 & 9769 & 9514 & 63  \\
			$15000\times15000$ && 14924 & 14869 & 14888 & 14968 & 68 \\
			\bottomrule
		\end{tabular}
		\label{tab1}
	\end{table}	
	
	In order to have a clearer view of the processing effect of the RCM algorithm on the structure of sparse matrices, we obtained a series of matrices from the Suite Sparse Matrix Collection (\url{https://sparse.tamu.edu/}). Taking matrices $poli3$, $blckhole$, $dixmaanl$, and $jagmesh4$ as examples, and we visualize the original structure of these matrices and the structure of these matrices after the reordering by the RCM algorithm, see Figure \ref{fig2}, where the blue part is the distribution of non-zero elements of the matrices. It can be clearly seen that after the RCM algorithm reordering, the nonzero elements of the matrices are closer to the diagonal, and the bandwidth of the matrices is significantly reduced.
	\begin{figure}[htbp]
		\centering
		\includegraphics[scale=0.45]{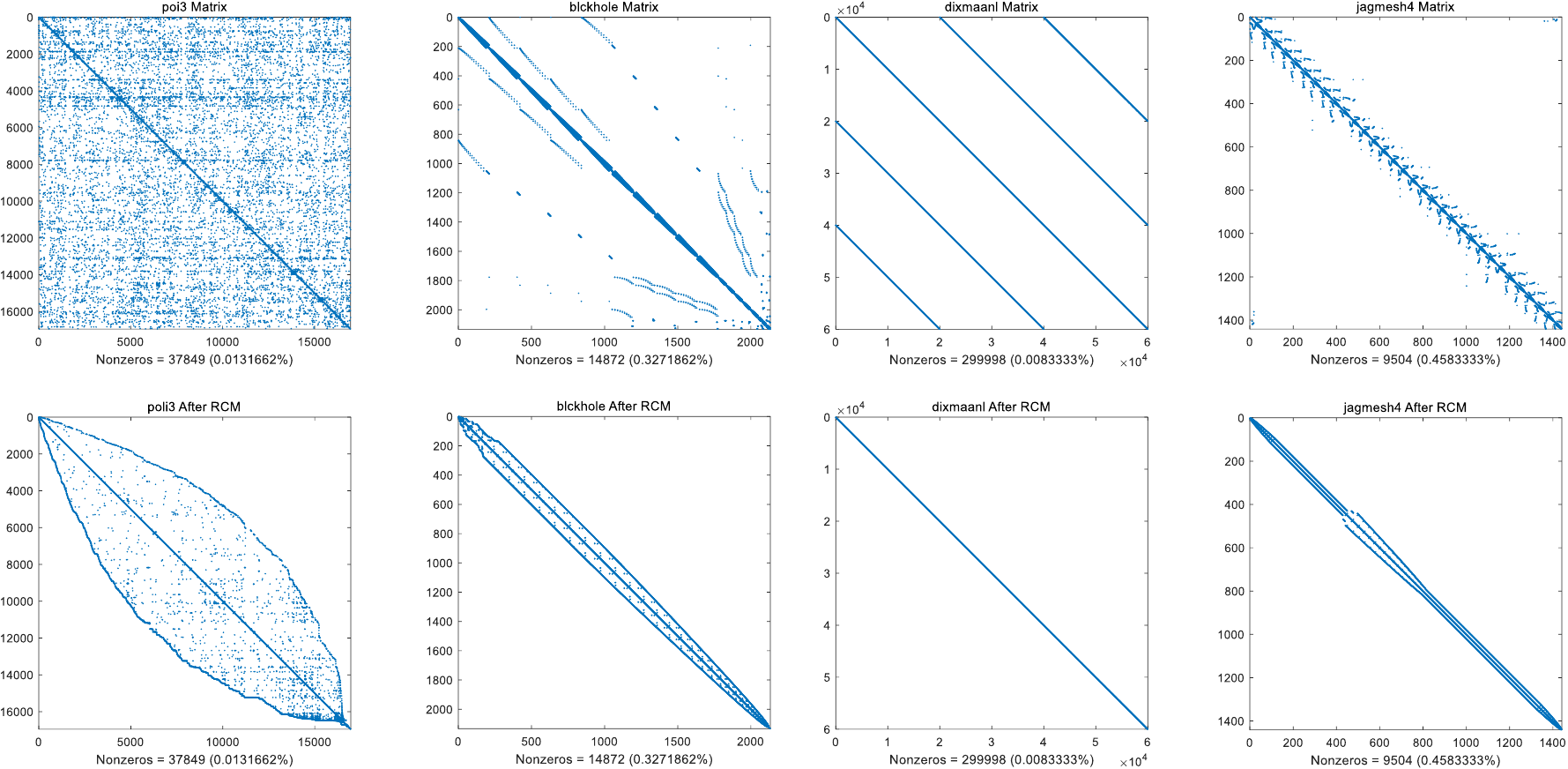}
		\caption{\; The top four images show the sparse structure of matrices poli3, blackhole, dixmaanl, and jagmesh4, while the bottom four images show the sparse structure of corresponding matrices after RCM Reordering.}
		\label{fig2}
	\end{figure}

	Meanwhile, the RCM algorithm is simple and efficient with low time complexity. That is, adding the RCM algorithm does not introduce significant time overhead. Taking the matrices $muu$, $poli3$, $rajat07$, $blckhole$, and $jagmesh4$ as an example, the time-consuming statistics of the RCM algorithm for reordering them are shown in the Table \ref{tab2}, where $bw$ denotes the original bandwidth of the matrices, and $bw1$ denotes the bandwidth of the matrices after reordering them by RCM. It can be seen that RCM has obvious advantages in reducing bandwidth and time consumption. For details on the time consumption of the RCM algorithm for reordering sparse matrices, see Ref. \cite{23}.
	\begin{table}[!hpt]
		\centering
		\setlength {\tabcolsep}{4mm}
		\setlength{\abovecaptionskip}{10pt}%
		\setlength{\belowcaptionskip}{10pt}%
		\caption{\;\;  Bandwidth variation and time consumption of the RCM algorithm for reordering sparse matrices.}
		\begin{tabular}{lcccccc}
			\toprule
			$Name$ && $m$ & $density\left(\%\right)$ & $bw$ & $bw1$ & $time~ consumption(s)$  \\
			\hline
			$muu$ && 7102 & 0.337 & 4696 & 311 & 0.0020  \\
			$poli3$ && 16955 & 0.013 & 16917 & 5487 & 0.0022  \\
			$rajat07$ && 14842 & 0.029 & 1483 & 413 & 0.0011  \\
			$blckhole$ && 2132 & 0.327 & 1805 & 105 & 0.0006  \\
			$jagmesh4$ && 1440 & 0.458 & 1408 & 56 & 0.0008  \\
			\bottomrule
		\end{tabular}
		\label{tab2}
	\end{table}	

	The above results show that the RCM sorting algorithm not only significantly reduces the bandwidth of the matrix, but also has a very low time overhead. Therefore, it is the preferred solution to directly use the RCM sorting algorithm to optimize the structure of sparse linear systems.
	
	\subsection{The Proposed methods}\label{sec3.3}
	As mentioned above, both RBK(k) and GREK(k) algorithms use the K-means algorithm to group matrix rows. Undoubtedly, K-means clustering is one of the most effective methods for dividing data into different groups. However, due to only considering the computational efficiency after partitioning and ignoring the computational cost of partitioning, the overall computational efficiency of the classic block Kaczmarz algorithm is not high. Therefore, in order to solve these problems, this paper proposes a simple and fast block Kaczmarz algorithm that combines matrix reordering techniques and orthogonal partitioning ideas.
	
	 Based on the above analysis, We reorder the original matrix $A$ as follows:
	\begin{equation}\label{eq:3.1}
		PAP^T,
	\end{equation} 
	where $P$ is obtained by the RCM reordering algorithm and $P$ is a substitution matrix (i.e., $P^TP=E$), so the linear system (1.1) becomes
	\begin{equation}\label{eq:3.2}
		PAP^TPx=Pf.
	\end{equation} 
	
	Let us denote $\widetilde{A}=PAP^T$, $\widetilde{f}=Pf$. Since  $\widetilde{A}$ is an approximate banded matrix, then it is sufficient to make the row blocks orthogonal to each other by doing appropriate row blocking of $\widetilde{A}$, see Figure \ref{fig2}.
	\begin{figure}[htbp]
		\centering
		\includegraphics[scale=0.5]{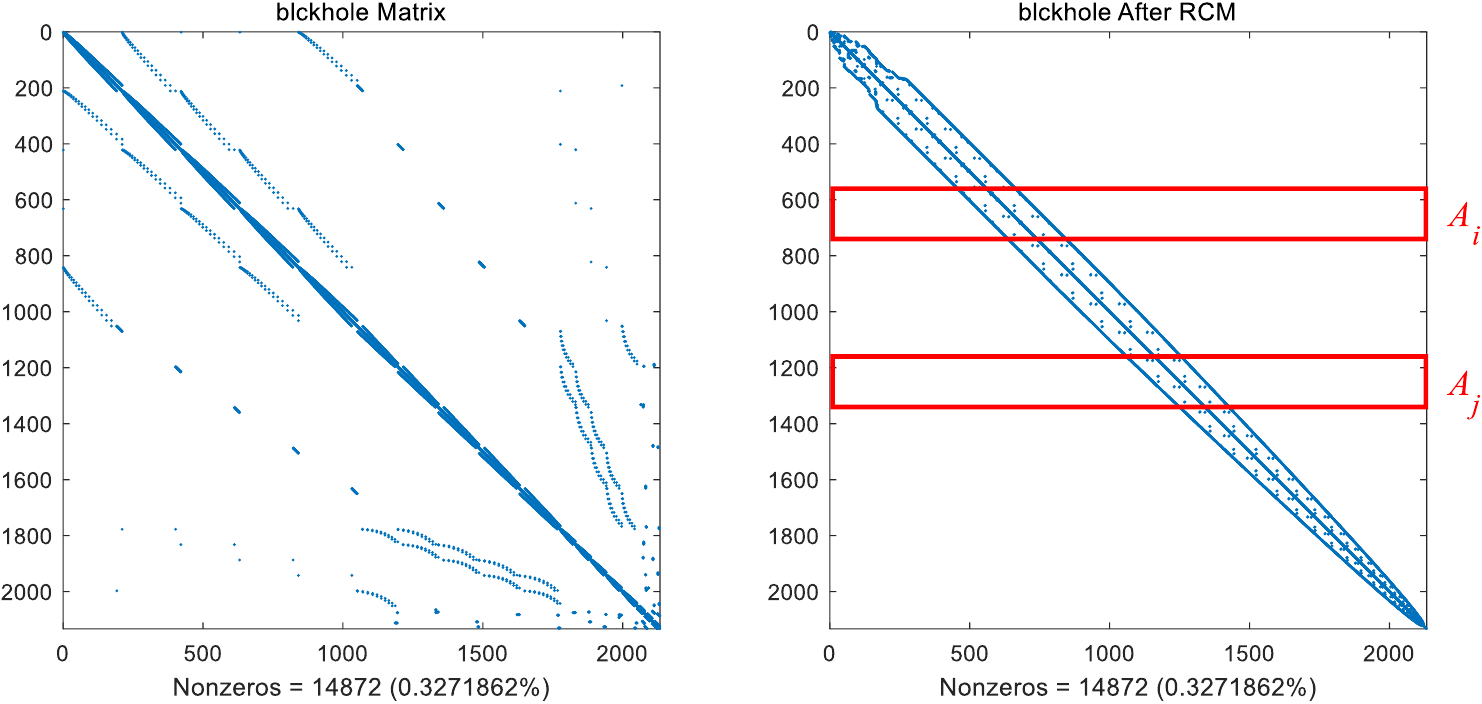}
		\caption{\; The orthogonality representation of the blckhole matrix after reordering with the RCM algorithm.}
		\label{fig2}
	\end{figure}

	Due to the natural orthogonality of banded matrices, we directly chunk $\widetilde{A}$ uniformly in order. Suppose that the matrix $\widetilde{A}$ is divided into $k$ blocks, the first $k-1$ blocks take $\left[m/k\right]$ elements each in order, and the remaining elements are stored in the $k$-th block, and the corresponding blocking is done for $\widetilde{f}$, 
	$$\widetilde{A}=\begin{bmatrix}
		A_1 \\
		A_2 \\
		...  \\
		A_k \\\end{bmatrix},~~~\widetilde{f}=\begin{bmatrix}
		f_1 \\
		f_2 \\
		...  \\
		f_k \\\end{bmatrix}.$$
	
   As is well known, the Kaczmarz method iteratively projects the current estimate orthogonally onto successive hyperplanes to approach the exact solution. As analyzed in Section 3.1, orthogonal hyperplanes enable accelerated convergence: pairwise orthogonal projections (e.g., $H_1\bot H_3$ in Figure \ref{fig1}) achieve quadratic convergence within two iterations, while non-orthogonal configurations require significantly more steps. This motivates a critical innovation: grouping mutually orthogonal hyperplanes into block projections to eliminate redundant iterations and accelerate convergence—a strategy central to our proposed method.

	Taking inspiration from this, we denote the centroid coordinates of each matrix block $A_i$ as $\overline{A}_i$, $i=1,2,...,k$, and then calculate the inner product $G\left(i,j\right)=\left(\overline{A}_i,\overline{A}_j\right)$, the two-paradigm number of $\overline{A}_i$ and $\overline{A}_j$, and then compute the cosine of the angle between the centroid coordinates according to the angle cosine formula $cos\theta_{ij}=\frac{G\left(i,j\right)}{\left\|\overline{A}_i\right\|_{2}\left\|\overline{A}_j\right\|_{2}}$, which constitutes a table $C$ of the cosine values, and $cos\theta_{ij}$ is used to measure orthogonality between blocks ${A}_i$ and blocks ${A}_j$. And since $\widetilde{A}$ is an approximate banded matrix, it is reasonable and convenient to use the cosine of the angle between the centroid coordinates to reflect the orthogonality between row blocks. And we divide the matrix $\widetilde{A}$ into $k$ blocks beforehand, and we only need to calculate the cosine of the angle between these $k$ blocks, and since $C$ is a symmetric matrix, we only need to calculate $\frac{k \left(k-1\right)}{2}$ times in practice. Since $k \ll m$, this is a fairly small number.
	
	Next, based on our analysis above, we aim to accelerate the Kaczmarz method by separating the two row blocks that are nearly orthogonal, i.e., $cos\theta_{ij}\approx 0$, into one class. Here we set a threshold $thr$ such that when $cos\theta_{ij}<thr$, block $A_{\tau_i}$ and block $A_{\tau_j}$ are considered orthogonal. However, it is obvious that for any block $A_{\tau_i}$, there is most likely more than one block that is orthogonal to it, and this can easily result in the situation where some rows blocks are in several classes at the same time. In order to avoid this situation, we stipulate that for each block $A_i$, $i=1,2,..k$, the first block $A_j \in \left\{A_i,A_{i+1},...,A_{k}\right\} $ orthogonal to $A_i$ is selected and placed in a class, noting the current class as $class\left\{i\right\}$, and requiring that there is no duplication of blocks between classes. We merge these orthogonal classes, denoted $Oclass$. For the blocks that do not go into $Oclass$, in order to avoid the non-convergence situation caused by the loss of information, we assign all of these blocks in a new class, denoted $Nclass$. This regulation is reasonable because of the nature of the A-approximate banded matrix. The solution is obtained by the block Kaczmarz method on both categories:
	\begin{equation}\label{eq:3.3}	\widetilde{x}_{j+1}=\widetilde{x}_j+A_{\tau_j}^{\dagger}\left(f_{\tau_j}-A_{\tau_j}\widetilde{x}_j\right),
	\end{equation} 
	performing two iterations on pairwise orthogonal blocks greatly accelerates the convergence speed of the Kaczmarz algorithm, and then performing Kaczmarz iterations on non-orthogonal blocks ensures that matrix information is not omitted.
	
	Finally, the solution of the original system is obtained as $x=P^T\widetilde{x}$. See Algorithm 4 for details.
	\begin{algorithm}[htb]
		\caption{ OBK-RCM Algorithm}
		\label{5}
		\begin{algorithmic}[1]  	
			\Require $A$, $f$, $x_0$, $k$, $thr$, $l$.
			\Ensure An estimation $x_{l}$ of the unique solution $x_{*}$ to $Ax = f$.
			\State Reordering matrices using the RCM algorithm, $\widetilde{A}=PAP^T$, $\widetilde{f}=Pf$.
			\State Perform uniform chunking of $\widetilde{A}$ and $\widetilde{f}$ in sequential order, and obtain $k$ groups $A_i$ and $f_i$, $i=1,2,...,k$. $\overline{A}=[\overline{A}_1,...,\overline{A}_k]$, $\overline{f}=[\overline{f}_1,...,\overline{f}_k]$, where $\overline{A}_i$ and $\overline{f}_i$ respectively represent the center of $A_i$ and $f_i$.
			\State Compute the cosine of the angle $cos\theta_{ij}$ between the centroid coordinates, and constitutes a table $C$ of the cosine values.
			\State Chunking $Oclass$ and $Nclass$ according to the algorithm description.
			\For {$j=0,1,2, \ldots$} until termination criterion is satisfied
			\State Select ~~$\tau_1,\tau_2 \in Oclass$
			\State Compute
			$$\widetilde{x}_{{j+1}/2}=\widetilde{x}_j+A_{\tau_1}^{\dagger}\left(f_{\tau_1}-A_{\tau_1}\widetilde{x}_j\right)$$
			$$\widetilde{x}_{j+1}=\widetilde{x}_{{j+1}/2}+A_{\tau_2}^{\dagger}\left(f_{\tau_2}-A_{\tau2}\widetilde{x}_{{j+1}/2}\right)$$
			\State Select ~~$\tau_{j} \in Nclass$
			\State Compute $$\widetilde{x}_{j+2}=\widetilde{x}_{j+1}+A_{\tau_j}^{\dagger}\left(f_{\tau_j}-A_{\tau_j}\widetilde{x}_{j+1}\right)$$
			\EndFor
			\State $x=P^T\widetilde{x}$.
		\end{algorithmic}
	\end{algorithm}
	
	Unlike the RBK, RBK(k) and GREBK(k) methods, the OBK-RCM algorithm does not require the construction of an integer set $U_j$ for ignoring row blocks corresponding to elements with small residuals. In fact, the calculation of coefficients $\epsilon_j$ and the construction of sets $U_j$ are very time-consuming in the code implementation process, which also makes the computational efficiency of the RBK (k) and GREBK (k) algorithms less ideal. Our proposed idea of performing $\frac{k \left(k-1\right)}{2}$ calculations of the angle cosine values and then orthogonal chunking avoids this problem to a large extent.On the one hand, the computation time of the angle cosine value is very short because $k \ll m$. On the other hand, two consecutive orthogonal projections onto mutually orthogonal hyperplanes also enable the iterative solution to quickly approach the true solution.
	
	\subsection{Convergence analysis}\label{sec3.4}
	Assuming that the exact solution of the original linear system $Ax=f$ is $x_*$, 
	and the exact solution of the reordered linear system $\widetilde{A}\widetilde{x}=\widetilde{f}$ is $\widetilde{x}_*$, where $P$ is the permutation matrix, then
	\begin{equation}\label{eq:3.4}
		\left\|x_j-x_*\right\|_{2}^{2}=\left\|P^T\widetilde{x}_j-P^T\widetilde{x}_*\right\|_{2}^{2}=\left\|P^T\left(\widetilde{x}_j-\widetilde{x}_*\right)\right\|_{2}^{2}=\left\|\widetilde{x}_j-\widetilde{x}_*\right\|_{2}^{2}.
	\end{equation} 
	It can be seen that the convergence analysis of the original linear system and the reordered linear system is consistent, therefore, we will not differentiate this in the future.	
	
	
	\textbf{Theorem 1.} If the linear system (\ref{eq:1.1}) $Ax=f$ is consistent, then the OBK-RCM algorithm produces a sequence of iterations $\left\{x_{j}\right\}_{j=0}^{\infty}$ that converges in expectation to a unique minimum norm solution $x_{\star}=A^{\dagger}f$. In particular, we have the following linear convergence rate in expectation:
	\begin{equation}\label{eq:3.5}
		\mathbb{E}\left\|x_j-x_*\right\|_{2}^{2} \leqslant \left(1-\frac{\lambda_{min}^{\tau}}{\lambda_{max}^{\tau}}\right)^{jk}\left\|x_0-x_*\right\|_{2}^{2}.
	\end{equation} 
	
	\textbf{Proof.} When $\tau_j \in Nclass$, the OBK-RCM algorithm involves one block per iteration and computes the linear system by the following form
	\begin{equation}\label{eq:3.6}
		{x}_{j+1}={x}_j+A_{\tau_j}^{\dagger}\left(f_{\tau_j}-A_{\tau_j}{x}_{j}\right).
	\end{equation} 
	According to the singular value decomposition, there are $A_{\tau_j}=U_{\tau_j}\Sigma_{\tau_j}V_{\tau_j}^H$, $A_{\tau_j}^{\dagger}=V_{\tau_j}\Sigma_{\tau_j}^HU_{\tau_j}^H$, and therefore $A_{\tau_j}^{\dagger}A_{\tau_j}=V_{\tau_j}MV_{\tau_j}^H$. Thus $A_{\tau_j}^{\dagger}A_{\tau_j}$ is an orthogonal projection operator, then
	\begin{equation}\label{eq:3.7}
		\begin{split}			x_{j+1}-x_*&=x_j+A_{\tau_j}^{\dagger}\left(f_{\tau_j}-A_{\tau_j}x_j\right)-x_* \\
		&= x_j-x_*+A_{\tau_j}^{\dagger}\left(A_{\tau_j}x_*-A_{\tau_j}x_j\right) \\ 
		&= \left(I-A_{\tau_j}^{\dagger}A_{\tau_j}\right)\left(x_j-x_*\right).
		\end{split}
	\end{equation}
	Therefore, we have
	\begin{equation}\label{eq:3.8}
		\begin{split}			\left\|x_{j+1}-x_*\right\|_{2}^{2}&=\left(x_j-x_*\right)^H\left(I-A_{\tau_j}^{\dagger}A_{\tau_j}\right)^H\left(I-A_{\tau_j}^{\dagger}A_{\tau_j}\right)\left(x_j-x_*\right) \\
			&=\left(x_j-x_*\right)^H\left(I-A_{\tau_j}^{\dagger}A_{\tau_j}\right)\left(x_j-x_*\right) \\
			&=\left\|x_j-x_*\right\|_{2}^{2}-\left\|A_{\tau_j}^{\dagger}A_{\tau_j}\left(x_j-x_*\right)\right\|_{2}^{2}.
		\end{split}
	\end{equation}
	Taking the conditional expectation on both sides of inequality (\ref{eq:3.8}), since the matrix $\widetilde{A}$ is uniformly partitioned, the probability of each block being selected is equal, i.e., $p_{\tau}=\frac{1}{k}$, we get
	\begin{equation}\label{eq:3.9}
		\begin{split}						\mathbb{E}_j\left[\left\|x_{j+1}-x_*\right\|_{2}^{2}\right]&= \left\|x_{j}-x_*\right\|_{2}^{2}-\mathbb{E}_j\left[\left\|A_{\tau_j}^{\dagger}A_{\tau_j}\left(x_j-x_*\right)\right\|_{2}^{2}\right] \\
			&= \left\|x_{j}-x_*\right\|_{2}^{2}-\sum_{\tau \in \left\{1,2,...,k\right\}}p_{\tau}\left[\left\|A_{\tau}^{\dagger}A_{\tau}\left(x_j-x_*\right)\right\|_{2}^{2}\right] \\
			&\leqslant \left\|x_{j}-x_*\right\|_{2}^{2}-\sigma_{min}^2\left(A_{\tau}^{\dagger}\right)\sigma_{min}^2\left(A_{\tau}\right)\left\|x_{j}-x_*\right\|_{2}^{2} \\
			&\leqslant \left(1-\frac{\sigma_{min}^2(A_{\tau})}{\sigma_{max}^2(A_{\tau})}\right)\left\|x_j-x_*\right\|_{2}^{2},
		\end{split}
	\end{equation}
	where $\sigma_{min}(.)$ denotes the minimum non-zero singular value of a matrix, $\sigma_{max}(.)$ denotes the maximum singular value of a matrix. It may be useful to set $\lambda_{min}^{\tau}$=$\sigma_{min}^2\left(A_{\tau}\right)$, $\lambda_{max}^{\tau}$=$\sigma_{max}^2\left(A_{\tau}\right)$. Then we have
	\begin{equation}\label{eq:3.10}
		\mathbb{E}_j\left[\left\|x_{j+1}-x_*\right\|_{2}^{2}\right] \leqslant \left(1-\frac{\lambda_{min}^{\tau}}{\lambda_{max}^{\tau}}\right)\left\|x_j-x_*\right\|_{2}^{2}.
	\end{equation}\label{eq:3.11}
	By taking the full expectation on both sides of the above equation, we have
	\begin{equation}
		\mathbb{E}\left\|x_{j+1}-x_*\right\|_{2}^{2} \leqslant \left(1-\frac{\lambda_{min}^{\tau}}{\lambda_{max}^{\tau}}\right)\mathbb{E}\left\|x_j-x_*\right\|_{2}^{2}.
	\end{equation}
	
	When $\tau_1, \tau_2 \in Oclass$, each iteration of the OBK-RCM algorithm involves two blocks, that is, iterates once on the basis of the original iteration, then obviously,
	\begin{equation}\label{eq:3.12}
		\begin{split}
			\mathbb{E}\left\|x_{j+1}-x_*\right\|_{2}^{2} &\leqslant \left(1-\frac{\lambda_{min}^{\tau}}{\lambda_{max}^{\tau}}\right)\mathbb{E}\left\|x_{(j+1)/2}-x_*\right\|_{2}^{2} \\
			&\leqslant \left(1-\frac{\lambda_{min}^{\tau}}{\lambda_{max}^{\tau}}\right)^2\mathbb{E}\left\|x_j-x_*\right\|_{2}^{2}.
		\end{split}
	\end{equation}
	
	We divide the matrix into $k$ blocks, where $k_1$ blocks are partitioned into $Oclass$ and $k_2$ blocks are partitioned into $Nclass$, and $k_1+k_2=k$. According to our POBK algorithm, iterating first in $Oclass$ and then in $Nclass$, then it is clear that the iteration is done once and 
	\begin{equation}\label{eq:3.13}
		\begin{split}
			\mathbb{E}\left\|x_{j+1}-x_*\right\|_{2}^{2} &\leqslant \left(1-\frac{\lambda_{min}^{\tau}}{\lambda_{max}^{\tau}}\right)^{k_1}\left(1-\frac{\lambda_{min}^{\tau}}{\lambda_{max}^{\tau}}\right)^{k_2}\mathbb{E}\left\|x_j-x_*\right\|_{2}^{2}  \\
			&\leqslant \left(1-\frac{\lambda_{min}^{\tau}}{\lambda_{max}^{\tau}}\right)^k\mathbb{E}\left\|x_j-x_*\right\|_{2}^{2}.
		\end{split}
	\end{equation}
	The inequality (\ref{eq:3.5}) follows by induction on $j$.
	\qed
	
	{\bf Remark 1.} From the above proof, it can be seen that the OBK-RCM algorithm converges as expected, and the rate of convergence is closely related to the maximum singular value of the block, the minimum non-zero singular value, and the number of chunks. When the submatrix orthogonality is enhanced (i.e., $\lambda_{max}^{\tau} \approx \lambda_{min}^{\tau}$), the ratio tends to 1 and the rate of convergence increases significantly. This explains the key role of orthogonal chunking of the matrix after reordering in accelerating the convergence.

	According to the analysis in Section \ref{sec3.1}, the influence of the angle between the matrix row blocks on the convergence speed is crucial. In order to visualize the convergence speed of the OBK-RCM algorithm more intuitively, the following section explains the reason why the OBK-RCM algorithm converges faster from the perspective of the angle between hyperplanes.
	
\begin{figure}[htbp]
		\centering
		\includegraphics[scale=0.5]{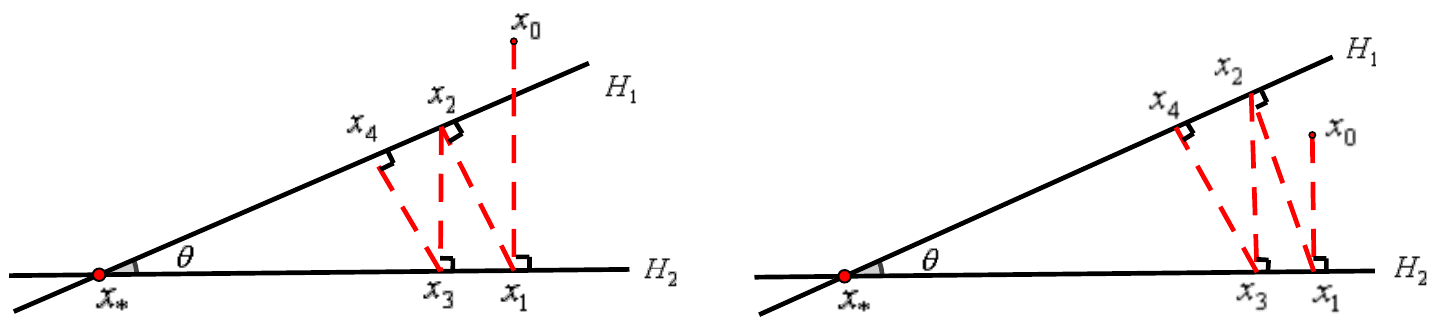}
		\caption{\;Two cases of the initial value $x_0$ on a 2D System.}
		\label{fig3}
	\end{figure}

	For a given linear system, given an initial value $x_0$, taking a two-dimensional linear system as an example, the Kaczmarz iteration method has two situations: the initial value $x_0$ is between two hyperplanes or outside of two hyperplanes, see Figure \ref{fig3}.

	Assuming the exact solution of the linear system is $x_*$, and the angle between hyperplanes $H_1$ and $H_2$ is $\theta$, for both of the above cases, we have that
	\begin{equation}\label{eq:3.14}
		cos\theta=\frac{\left\|{x_2-x_*}\right\|_{2}}{\left\|{x_1-x_*}\right\|_{2}}=\frac{\left\|{x_3-x_*}\right\|_{2}}{\left\|{x_2-x_*}\right\|_{2}}=...=\frac{\left\|{x_j-x_*}\right\|_{2}}{\left\|{x_{j-1}-x_*}\right\|_{2}},
	\end{equation}
	thus
	\begin{equation}\label{eq:3.15}
		\left\|x_j-x_*\right\|_{2}=\left\|x_{j-1}-x_*\right\|_{2}cos\theta=...=\left\|x_1-x_*\right\|_{2}cos^{j-1}\theta.
	\end{equation}
	
	According to Theorem 1, so
	\begin{equation}\label{eq:3.16}
		\left\|x_j-x_*\right\|_{2}^{2}\leq\left(1-\frac{\lambda_{min}^{\tau}}{\lambda_{max}^{\tau}}\right) \cdot cos^{2j-2}\theta \cdot \left\|x_0-x_*\right\|_{2}^{2}.
	\end{equation}
	
	From Equation (\ref{eq:3.16}), it can be seen that the convergence speed of the Kaczmarz algorithm is closely related to the cosine of the hyperplane angle. The larger the hyperplane angle $\theta$ is (the closer to $90^{\circ}$), and the closer $cos\theta$ is to 0, the faster the algorithm converges; on the other hand, the smaller the hyperplane angle  $\theta$ is (the closer to $0^{\circ}$), and the closer $cos\theta$ is to 1, the slower the convergence speed of the algorithm is. The OBK-RCM algorithm utilizes the RCM sorting algorithm to reorder the matrices so that the hyperplanes $H_1$ and $H_2$ are approximately orthogonal ($cos\theta \approx 0$), and the error is rapidly decayed so that the convergence speed is rapidly improved.
	
	Briefly, the OBK-RCM algorithm achieves superior convergence over traditional block Kaczmarz methods (RBK, RBK(k), GREBK(k), aRBK) through two synergistic mechanisms:
	
	\begin{enumerate}
	  \item RCM-Driven Orthogonality Enhancement: By reordering the matrix into a banded structure (Table 1), OBK-RCM forces inter-block angles toward orthogonality, enabling exact solution convergence in two iterations for orthogonal hyperplanes (Figure \ref{fig1}), versus polynomial-time scaling for non-orthogonal cases.
	  \item Intelligent Block Scheduling: The $Oclass$/$Nclass$ partitioning (Section 3.3) systematically minimizes adjacency of non-orthogonal blocks in $Nclass$. Although $Nclass$ blocks lack strict orthogonality, their residual angular cosines remain significantly lower than in the original matrix $A$ (Equation \ref{eq:3.16}), as RCM reordering inherently suppresses small-angle configurations (Figure 2).
	\end{enumerate}
	
	This dual strategy--orthogonal acceleration in $Oclass$ and mitigated error propagation in $Nclass$--collectively ensures OBK-RCM’s 10-50 times speedup over baseline methods, as quantified in Table 4.

	\subsection{Numerical experiments}\label{sec3.5}
	In this section, we conducted some experiments to compare the convergence speed of RBK, RBK(k), GREBK(k), aRBK and OBK-RCM methods, verifying the effectiveness and efficiency of the algorithm. The matrices used for numerical experimental tests are 15 sparse matrices selected from the $SuiteSparse \; Matrix \;$ $ Collection$ with application background and varying in size, condition number and density. Note that when matrix $A$ is a low rank or underdetermined matrix, the Matlab function pinv often cannot obtain an accurate solution to the linear system. We use the Matlab function lsqminnorm to solve the linear equation $Ax=f$ and minimize the value of norm in the vector space.
	
	For all the above algorithms, we set the initial value $x_0=0$ and the exact solution is $x_*$. If the relative error $RSE$ under the current iteration is less than $1e-6$ or the number of iterations exceeds $5e+5$, the calculation will be stopped, where
	$$RSE=\frac{\left\|x_j-x_*\right\|_{2}^2}{\left\|x_*\right\|_{2}^2}.$$
	
	All experiments were conducted in $MATLAB (R2019a)$ on a PC with AMD Ryzen 7 5825U with Radeon Graphics 2.00 GHz. We use iteration steps (IT) and CPU time in seconds (CPU) to evaluate the numerical performance of different algorithms. Here, IT and CPU are the average iteration steps and CPU time of running the underlying algorithm 10 times. If the corresponding algorithm fails to converge under the specified relative error bound or maximum number of iterations, the corresponding IT and CPU are recorded as Inf and NAN in the table.
	
	We also provided the basic information of the test matrix, including matrix dimension $\left(m \times m\right)$, matrix condition number, and matrix density, as detailed in Table \ref{tab3}. The definition of matrix density is:
	$$density=\frac{number~of~nonzeros~of~an~m \times m~matrix}{m \cdot m}.$$
	\begin{table}[hpt]
		\centering
		\setlength {\tabcolsep}{1.5mm}
		\setlength{\abovecaptionskip}{5pt}%
		\setlength{\belowcaptionskip}{10pt}%
		\caption{\;Information on the test matrices from the $SuiteSparse \; Matrix \; Collection$.}	
		\begin{tabular}{lclrrcc}
			\toprule
			$Name$ && $Kind$ & $m$ & $Nonzeros$ & $Cond$ & $ density\left(\%\right)$ \\
			\hline
			$poli$ && $Economic~Problem$ & 4008 & 8188 & 3.115e+02 & 0.051 \\
			$muu$ && $Structural~Problem$ & 7012 & 170134 & 7.654e+01 & 0.337 \\
			$poli3$ && $Economic~Problem$ & 16955 & 37849 & 5.253e+02 & 0.013 \\
			$ex29$ && $Computational~Fluid~Dynamics~Problem$ & 2870 & 23754 & 7.515e+02 & 0.288 \\
			$kim1$ && $2D/3D~Problem$ & 38415 & 933195 & 9.863e+03 & 0.063 \\
			$qpband$ && $Optimization~Problem$ & 20000 & 45000 & 6.436e+00 & 0.011 \\
			$rajat07$ && $Circuit~Simulation~Problem$ & 14842 & 63913 & 7.894e+02 & 0.029 \\
			$blckhole$ && $Structural~Problem$ & 2132 & 14872 & 4.167e+03 & 0.327 \\
			$linverse$ && $Statistical/Mathematical~Problem$ & 11999 & 95977 & 3.947e+03 & 0.067 \\
			$torsion1$ && $Duplicate~Optimization~Problem$ & 40000 & 197608 & 4.099e+01 & 0.012 \\
			$polilarge$ && $Economic~Problem$ & 15575 & 33074 & 2.983e+01 & 0.014 \\
			$jagmesh4$ && $2D/3D~Problem$ & 1440 & 9504 & 1.516e+04 & 0.458 \\
			$bcsstm39$ && $Structural~Problem$ & 46772 & 46772 & 8.271e+03 & 0.002 \\
			$crystm03$ && $Materials~Problem$ & 24696 & 583770 & 2.640e+02 & 0.096 \\
			$chem97ztz$ && $Statistical/Mathematical~Problem$ & 2541 & 7361 & 2.472e+02 & 0.114 \\
			\bottomrule
		\end{tabular}
		\label{tab3}
	\end{table}	

	Table \ref{tab4} gives the number of iterations and computation time for these test matrices under different algorithms. Since the RBK, GREBK(k) and RBK(k) algorithms involve only one block of data per iteration, whereas the aRBK and OBK-RCM algorithms involve all the blocks of data of the matrix per iteration, it is not fair to directly compare the number of iteration steps (IT) of the different algorithms. As can be seen from table \ref{tab4}, even if the OBK-RCM algorithm multiplies the number of iteration steps by $k$ times ($k<20$), it still outperforms the RBK, GREBK(K), and RBK(k) algorithms. In terms of computation time, the OBK-RCM algorithm is much faster than the randomized block algorithm Kaczmarz (RBK, aRBK) and the K-means direct block algorithm Kaczmarz(GREBK(k), RBK(k)). It is worth noting that for the matrix kim1 in the complex domain, the RBK, GREBK(k), and RBK(k) algorithms encountered problems during the chunking process and therefore could not continue to run. However, the OBK-RCM algorithm runs smoothly and achieves the desired convergence rate. This shows that it is not only applicable to the real domain but also to the complex domain.
	\begin{table}[!hpt]
		\centering
		\setlength {\tabcolsep}{4mm}
		\setlength{\abovecaptionskip}{5pt}%
		\setlength{\belowcaptionskip}{10pt}%
		\caption{\; IT and CPU performance of 15 test matrices under different algorithms}
		\resizebox{\textwidth}{102mm}{
			\begin{tabular}{llllllll}
				\toprule
				$Name$ &\ & $RBK$ & $RBK(k)$ & $GREBK(k)$& $aRBK$ & $OBK-RCM$ \\
				\hline
				$poli$ & IT & 74566 & 38636 & 51259 & 1721 & 28  \\
				& CPU & 40.63 & 43.43 & 28.71 & 33.75 & 0.07  \\
				\midrule
				$muu$ & IT & 160 & 53 & 269 & 13 & 1 \\
				& CPU & 8.44 & 5.38 & 2.44 & 3.51 & 0.22  \\
				\midrule
				$poli3$ & IT & 26621 & 6615 & 18902 & 1140 & 1034 \\
				& CPU & 104.62 & 125.98 & 48.07 & 103.04 & 7.38  \\
				\midrule
				$ex29$ & IT & 123 & 10 & 116 & 7 & 8 \\
				& CPU & 0.65 & 1.37 & 0.20 & 0.51 & 0.10  \\
				\midrule
				$kim1$ & IT & NAN & NAN & NAN & 118 & 52 \\
				& CPU & NAN & NAN & NAN & 276.98 & 18.27  \\
				\midrule
				$qpband$ & IT& 1687 & 43 & 292 & 11 & 1 \\
				& CPU & 3.31 & 12.33 & 4.33 & 3.36 & 1.01  \\
				\midrule
				$rajat07$ & IT & Inf & 38150 & 42957 & 3478 & 1080 \\
				& CPU & NAN & 2266.82 & 1889.06 & 585.08 & 57.52  \\
				\midrule
				$blckhole$ & IT & 48327 & 22213 & 198250 & 2308 & 1055 \\
				& CPU & 194.91 & 32.55 & 188.74 & 151.08 & 10.10  \\
				\midrule
				$linverse$ & IT & Inf & Inf & Inf & Inf & 2 \\
				& CPU & NAN & NAN & NAN & NAN & 0.42  \\
				\midrule
				$torsion1$ & IT & 37424 & 533 & 845 & 386 & 132 \\
				& CPU & 166.10 & 119.36 & 32.65 & 209.02 & 18.57  \\
				\midrule
				$polilarge$ & IT & 125195 & 109525 & 70322 & 11102 & 5128 \\
				& CPU & 411.08 & 257.46 & 149.94 & 625.90 & 20.40  \\
				\midrule
				$jagmesh4$ & IT & 385621 & 125814 & 73429 & 7104 & 2873 \\
				& CPU & 1173.66 & 130.63 & 96.18 & 292.11 & 10.40  \\
				\midrule
				$bcsstm39$ & IT & 50 & 2 & 10 & 4 & 1 \\
				& CPU & 0.16 & 122.52 & 9.30 & 4.09 & 0.04  \\
				\midrule
				$crystm03$ & IT & 4503 & 21 & 196 & 29 & 1 \\
				& CPU & 41.89 & 36.45 & 16.20 & 29.81 & 2.15  \\
				\midrule
				$chem97ztz$ & IT & 8727 & 2037 & 3574 & 20 & 9 \\
				& CPU & 3.88 & 3.64 & 2.03 & 0.61 & 0.04  \\
				\bottomrule
		\end{tabular}}
		\label{tab4}
	\end{table}

\begin{figure}[htbp]
		\centering
		\includegraphics[scale=0.45]{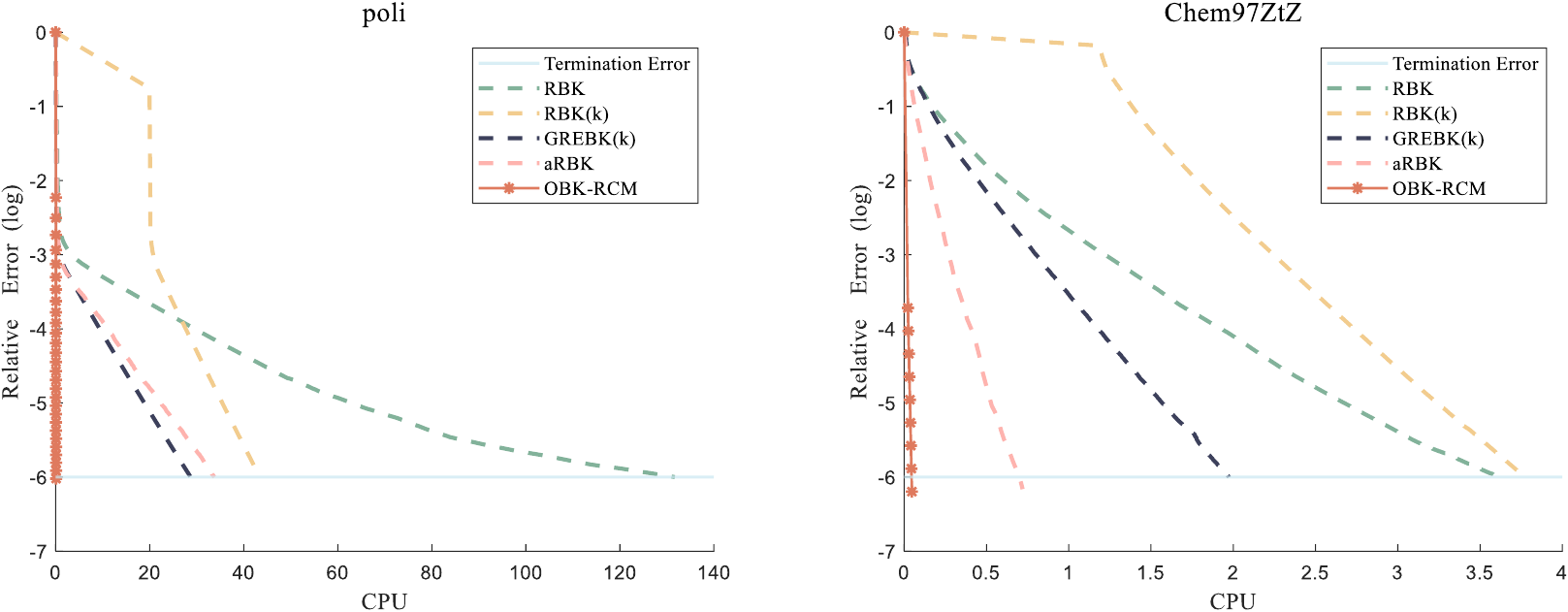}
		\caption{\;Average computational time for test matrices $poli$ and $torsion1$ to reach convergence for different iteration methods.}
		\label{fig4}
	\end{figure}

	From the Table \ref{tab4} and Figure \ref{fig4}, it can be seen that the OBK-RCM algorithm is able to complete the computation within a smaller number of iterations and a shorter CPU time, which exhibits extremely high computational efficiency and stability. This indicates that the OBK-RCM algorithm has significant advantages in dealing with large-scale linear systems with sparse square matrices, especially when the structure of the sparse matrices is more decentralized.
	
	Theoretically, the OBK-RCM algorithm has a better convergence rate by placing $x=P^T\widetilde{x}$ inside the loop, which implies that the OBK-RCM algorithm still has potential.
	
	Experimental results (Table \ref{tab4}, Figures \ref{fig4}) demonstrate OBK-RCM's dual superiority:
	\begin{itemize}
	  \item Accelerated Convergence: Requires 50--90\% fewer iterations than RBK/GREBK(k) across matrices with dispersed nonzero patterns (e.g., blckhole: 1,055 vs 198,250 iterations).
	  \item Runtime Efficiency: Achieves 10--50 times CPU time (up to several hundred) reduction by eliminating costly K-means clustering and residual sampling (Section 3.3).
	\end{itemize}

	This performance stems from RCM's banded structure optimization and Oclass/Nclass dynamic partitioning—strategies that jointly maximize hyperplane orthogonality while minimizing small-angle configurations, as quantified by Equation (\ref{eq:3.16}).

	\subsection{Select the optimal number of blocks of the matrix}\label{sec3.7}
	As mentioned earlier, for sparse linear systems, blocking is one of the effective means to improve the convergence of kaczmarz algorithm. Taking the $crystm03$ matrix as an example, Figure \ref{fig6} shows IT and CPU required for the matrix to reach convergence for different values of $k$. It can be seen that when $k=3$, the OBK-RCM algorithm has the optimal number of iterations and running time. 
	\begin{figure}[htbp]
		\centering
		\includegraphics[scale=0.45]{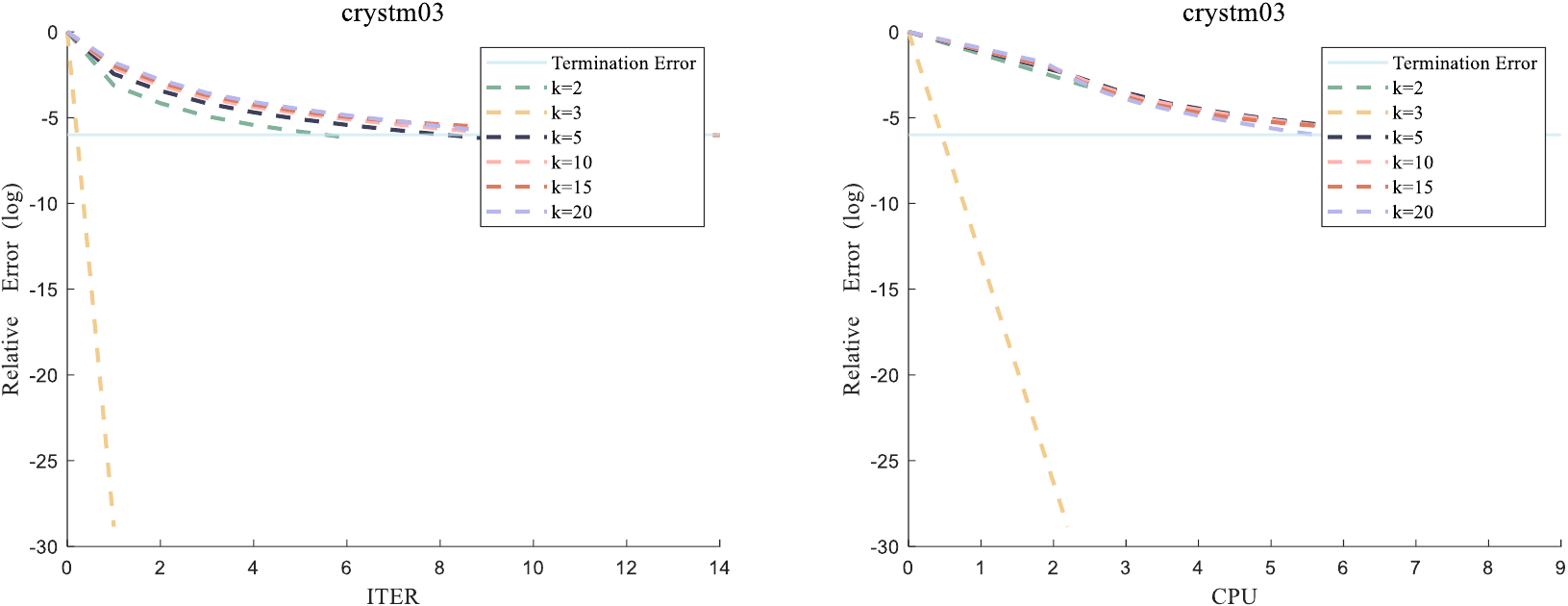}
		\caption{\;IT and CPU required for matrix $crystm03$ to reach convergence for different values of $k$.}
		\label{fig5}
	\end{figure}

	If $k$ is too large, the dimension of each block becomes smaller, and although the projection computation cost is reduced, the orthogonality between rows of blocks is reduced (i.e., the angle between blocks increases), which slows down convergence; if $k$ is too small, the dimension of a single block is too large, the computational complexity within the block rises and the accelerating effect of the hyperplane orthogonality on the speed of convergence cannot be fully utilized. A proper $k$ maximizes the minimum singular value and thus increases the convergence speed.
	
	According to Equation (\ref{eq:3.16}), the convergence speed of the OBK-RCM algorithm is directly related to the angle between rows of blocks, where the angle $\theta$ involves both mutually orthogonal and non-mutually orthogonal blocks. When mutually orthogonal blocks are grouped in a category, the cosine of the angle between these mutually orthogonal blocks will be very small or even close to 0. At this point, it is important to focus on the non-mutually orthogonal blocks, i.e., the number of nonzero elements in the cosine table, and the size of these nonzero elements for the rate of convergence of the OBK-RCM algorithm to be crucial to the number of matrix blocks. Here, two variables are defined: the orthogonality ratio $zn$ and the degree of non-orthogonality $nn$ to be used as selection criteria for the number $k$ of chunks.
	
	\textbf{Definition  1.} For the cosine value table $C$,  if the total number of elements is $num=k \times k$, the number of zero elements is $n_1$, then the Orthogonality Proportion $zn$ is denoted by
	$$zn=\frac{n_1}{num}=\frac{n_1}{k \times k}.$$
	
	\textbf{Definition  2.} For the cosine value table $C$, the total number of elements is $num=k \times k$, the number of non-zero elements is $n_2$, and the average of these non-zero elements is $nm$, then the Non-orthogonality Severity $nn$ is denoted by
	$$nn=\frac{n_2 \times nm}{num}=\frac{n_2 \times nm}{k \times k}.$$
	Theoretically, the optimal $k$ should simultaneously maximize $zn$ and minimize $nn$.
	
	It is worth noting that the variable $nn$ is not directly defined by the number of non-zero elements in the matrix $C$ as a percentage of the total number of elements, because we have to consider not only the number of non-zero elements, but also the size of the values of these non-zero elements. If these non-zero elements have large values, it means that the angle between some blocks is small, which also makes the iteration slower, which is exactly what we do not want to see.
	
	We take the $muu$ matrices as examples. 
	For the $muu$ matrix, from Figure  \ref{fig6}, $zn$ obtains its maximum value at $k=2$, and it happens that $nn$ also obtains its minimum value at $k=2$, so the choice of $k = 2$ is appropriate, and the lowest point of the running time box-and-line plot is also at $k = 2$ from Figure \ref{fig7}.
	\begin{figure}[htbp]
		\centering
		\includegraphics[scale=0.48]{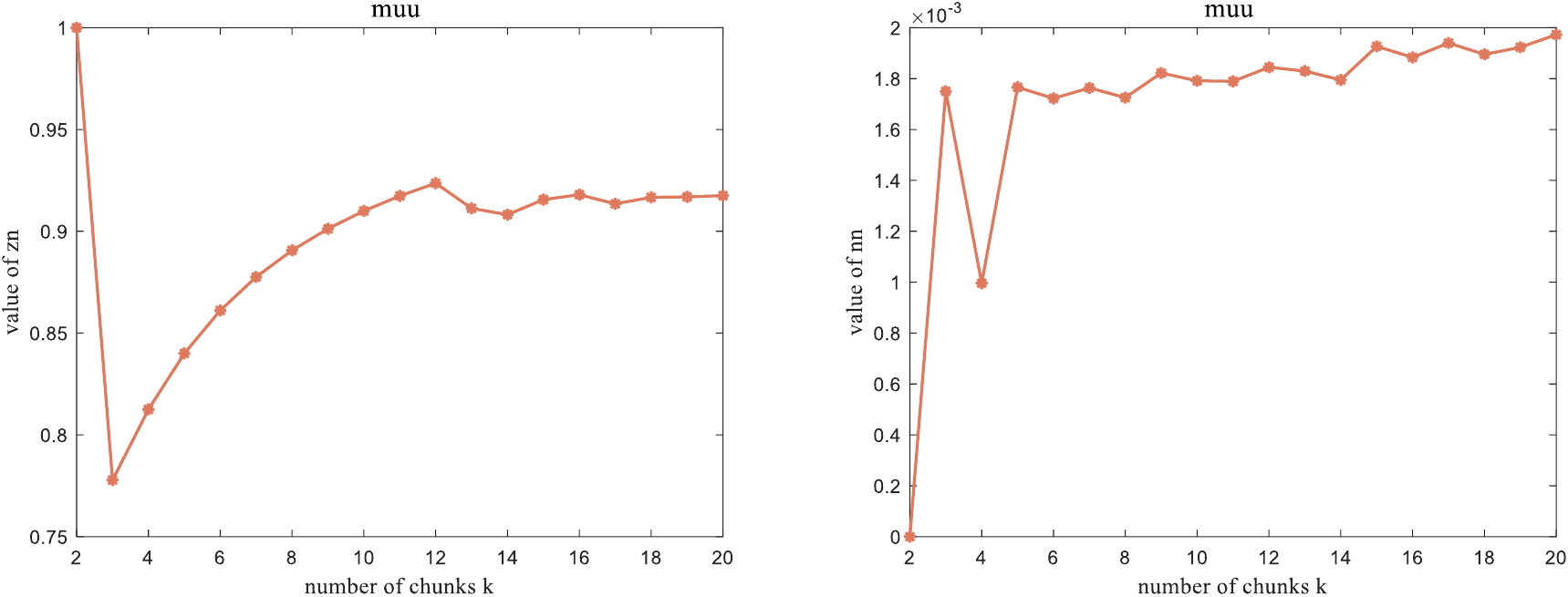}
		\caption{\; $zn$ and $nn$ values of the $muu$ matrix for different values of $k$.}
		\label{fig6}
	\end{figure}
	\begin{figure}[htbp]
		\centering
		\includegraphics[scale=0.55]{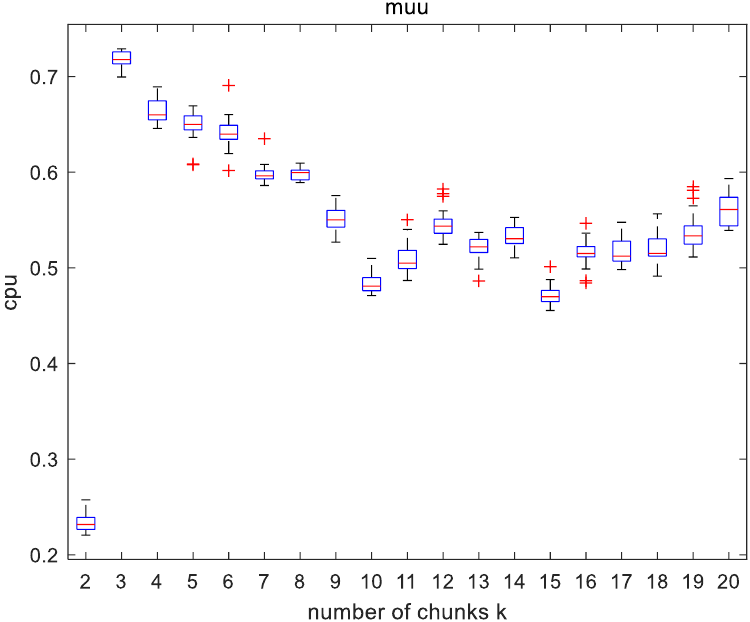}
		\caption{\;Boxplots of runtime for $muu$ matrices for different values of $k$.}
		\label{fig7}
	\end{figure}

	\subsection{Determination of threshold $thr$}\label{sec3.7}
	According to the previous definition, the threshold $thr$ is a key parameter in determining the orthogonality of matrix rows and blocks. Too large or too small a threshold value will reduce the convergence speed of the OBK-RCM algorithm. Therefore, a reasonable choice of the threshold $thr$ requires a trade-off between convergence speed and computational efficiency.
	
	In this experiment, in order to balance the computational time and computational effectiveness, the running times of typical matrices ($ex29$, $blckhole$, $jagmesh4$ and $linverse$) with different $thr$ values are compared, and the $thr$ value with optimal performance for most of the matrices is selected as the global threshold setting. Setting $thr1$, $thr2$, $thr3$, $thr4$, $thr5$, and $thr6$ as $0.5$, $0.2$, $0.1$, $0.05$, $0.02$, and $0.01$, respectively, it can be seen from Table \ref{tab5} that setting $thr5$ (i.e., setting $thr = 0.02$) as the global threshold is worthwhile.
	\begin{table}[!hpt]
		\centering
		\setlength {\tabcolsep}{4mm}
		\setlength{\abovecaptionskip}{10pt}%
		\setlength{\belowcaptionskip}{10pt}%
		\caption{\;\;Comparison of processing matrix bandwidth of different sorting algorithms.}
		\begin{tabular}{ccccccccc}
			\toprule
			$Name$ && $thr1$ & $thr2$ & $thr3$ & $thr4$ & $thr5$ & $thr6$  \\
			\hline
			$ex29$ && 0.1126 & 0.1097 & 0.1087 & 0.1086 & 0.1012 & 0.1014 \\
			$linverse$ && 0.4295 & 0.4315 & 0.4281 & 0.4259 & 0.4265 & 0.4281  \\
			$blckhole$ && 10.4855 & 10.257 & 10.1727 & 10.2027 & 10.1444 & 10.2011  \\
			$jagmesh4$ && 10.5939 & 10.6184 & 10.5471 & 10.5426 & 6810.5386 & 10.5355 \\
			\bottomrule
		\end{tabular}
		\label{tab5}
	\end{table}	

	\subsection{Stability of the OBK-RCM algorithm}\label{sec3.8}
	The K-means algorithm, as one of the most basic clustering methods, has the advantages of simplicity and high efficiency, but the results of this algorithm vary each time the clustering is performed, which is reflected in the numerical experiments as the block Kaczmarz algorithm for K-means clustering, where the computation time is very different even for the same $k$ value, which leads to the algorithm's stability is not very good. 
	\begin{figure}[htbp]
		\centering
		\includegraphics[scale=0.55]{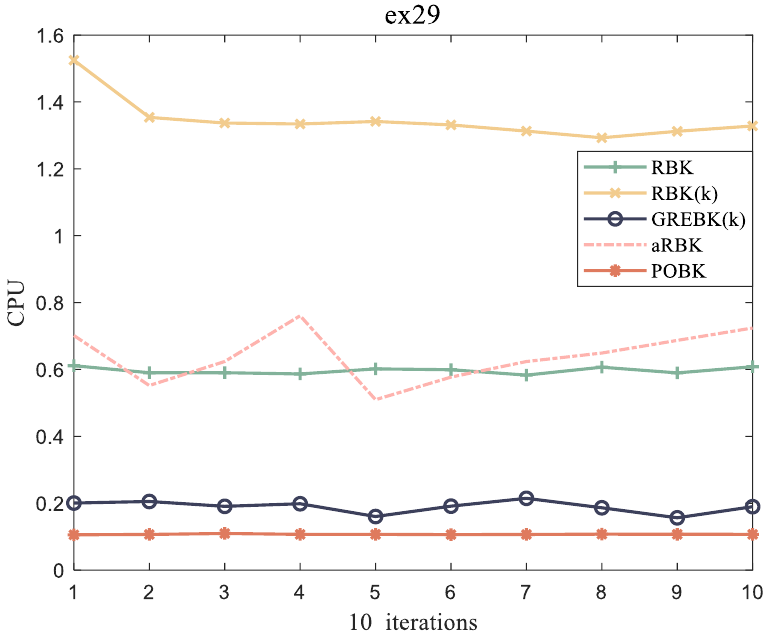}
		\caption{\;Changes in iteration time for 10 iterations of different algorithms under constant $k$.}
		\label{fig8}
	\end{figure}
	
	In the method proposed in this section (OBK-RCM), given a test matrix $A$, the approximate banded matrix $A$ after reordering by the RCM algorithm is unique, and since the OBK-RCM algorithm performs sequential uniform chunking, the form of the chunks is fixed once $k$ is determined, and hence the algorithm is stable for the same value of $k$. Numerical experiments also show that the algorithm is stable, and taking the ex29 matrix as an example, Figure \ref{fig8} shows the iteration time variations of the RBK, RBK(k), GREBK(k), aRBK, and OBK-RCM algorithms after 10 iterations with $k$ constant. It can be clearly seen that the OBK-RCM algorithm not only outperforms the other three algorithms in terms of computation time, but also performs better in terms of computation time stability.

	\section{Extended Orthogonal Block Kaczmarz for Non-Square Systems}
  The OBK-RCM algorithm, introduced in Section \ref{sec3}, achieves efficient solutions for sparse linear systems with scattered nonzero patterns via reordering and orthogonal block partitioning. However, a key limitation remains: OBK-RCM requires the coefficient matrix $A$ to be square $\left(i.e., m=n\right)$. To address this, we propose extended matrix formulations that preserve sparsity while extending compatibility to rectangular systems. These extensions retain RCM’s orthogonality-enhancing properties and $Oclass/Nclass$ scheduling, ensuring no computational efficiency loss (Figures \ref{fig9}-\ref{fig11}).

	Next, the compatible sparse non-square linear system $Ax = f$ is considered, where $A$ is an $m \times n$ matrix $\left(m \neq n\right)$. This section will analyze and discuss both $m > n$ and $m < n$.
	
	\subsection{The $m>n$ case}\label{sec4.1}
	For the linear system $Ax = f$, when $m > n$, OBK-RCM fails because the input needs to be square. For this reason, a new matrix $\widehat{A}$ is constructed as:
	\begin{equation}
	\begin{aligned}\label{eq:4.11}
		\widehat{A}=[A,O]
	\end{aligned}
	\end{equation}
	where $O$ is the $m \times (m - n)$ zero matrix. Thus $\widehat{A}$ becomes an $m \times m$ square matrix.  
	
	Also, let $\widehat{f}=f$, for the equation $\widehat{A}\widehat{x}=\widehat{f}$:
	\begin{equation}
	\begin{aligned}\label{eq:4.11}
				\widehat{A}\widehat{x}
		= \left[A, O\right] \cdot
		\begin{bmatrix}
			\widehat{x}_1 \\
			\widehat{x}_2 \\
		\end{bmatrix}
		=A\widehat{x}_1+0'
		=f
		=\widehat{f}	
	\end{aligned}
	\end{equation}
	where $\widehat{x}_1$ is an $n\times 1$ vector, $\widehat{x}_2$ is an $(m-n)\times 1$ vector, $0'$ is an $m\times 1$ zero vector.
	
	Since the linear system $Ax = f$ is compatible, it follows that $\widehat{A}\widehat{x}=\widehat{f}$ is also compatible. In fact, the solution of $\widehat{A}\widehat{x}=\widehat{f}$ has a larger solution set than $Ax = f$ because $\widehat{x}_2$ can take on any value. 
	
	Based on the above analysis, it can be shown that the first $n$ rows of the solution $\widehat{A}\widehat{x}=\widehat{f}$ is the solution of the linear system $Ax = f$. The specific algorithm is described in Algorithm 5.
	\begin{algorithm}[htb]
		\caption{\;\; SOBK-RCM Algorithm}
		\label{alg 5}
		\begin{algorithmic}[1]  	
			\Require The sparse matrix $A$, a right-hand term $f$, the initial value $x_{0}$, the number of blocks $k$, thresholds $thr$, the maximum number of iterations $l$.
			\Ensure An estimation $x_{l}$ of the unique solution $x_{*}$ to $Ax = f$.
			\State Construct a new matrix $\widehat{A}=[A,O]$, and let $\widehat{f}=f$. 
			\State For the new linear system $\widehat{A}\widehat{x}=\widehat{f}$, the OBK-RCM algorithm is used to compute the solution $\widehat{x}$.
			\State $x=\widehat{x}(1:n)$ is the solution of the original linear system $Ax = f$, where $n$ is the dimension of the original matrix $A$. 
		\end{algorithmic}
	\end{algorithm}
	
	
	\subsection{The $m<n$ case}\label{sec4.2}
	For the linear system $Ax = f$, when $m < n$, OBK-RCM fails because the input needs to be square. For this reason, a new matrix $\widehat{A}$ is constructed as:
	\begin{equation}
		\widehat{A}=\begin{bmatrix}
			A \\
			O \\
		\end{bmatrix}
	\end{equation}
	where $O$ is the $(n-m)\times n$ zero matrix. Thus $\widehat{A}$ becomes an $n \times n$ square matrix.
	
	Also, let 	$\widehat{f}=\begin{bmatrix}
		f \\
		0' \\
	\end{bmatrix}$, where $0'$ is $(n-m)\times 1$ zero vector , for the equation $\widehat{A}\widehat{x}=\widehat{f}$:
	\begin{equation}
    	\widehat{A}\widehat{x}
		= \begin{bmatrix}
			A \\
			O \\
		  \end{bmatrix} \cdot \widehat{x}
		=\begin{bmatrix}
			A\widehat{x} \\
			0' \\
	    	\end{bmatrix}
		=\begin{bmatrix}
			f \\
			0' \\
		  \end{bmatrix}
		=\widehat{f}.	
	\end{equation}

	Obviously, the solution $\widehat{x}$ of the extended systgem $\widehat{A}\widehat{x}=\widehat{f}$ satisfies $\widehat{x}=x_*$, where $x_*$ is the exact solution of $Ax = f$. The specific algorithm is described in Algorithm 6.
	\begin{algorithm}[htb]
		\caption{ UOBK-RCM Algorithm}
		\label{5}
		\begin{algorithmic}[1]  	
			\Require The sparse matrix $A$, a right-hand term $f$, the initial value $x_{0}$, the number of blocks $k$, thresholds $thr$, the maximum number of iterations $l$.
			\Ensure An estimation $x_{l}$ of the unique solution $x_{*}$ to $Ax = f$.
			\State Construct a new matrix $\widehat{A}=[A;O]$, and let $\widehat{f}=[f;0']$. 
			\State For the new linear system $\widehat{A}\widehat{x}=\widehat{f}$, the OBK-RCM algorithm is used to compute the solution $\widehat{x}$.
			\State $x=\widehat{x}=P^T\widetilde{x}$ is the solution of the original linear system $Ax = f$. 
		\end{algorithmic}
	\end{algorithm}
	
	With the above analysis, the SOBK-RCM and UOBK-RCM algorithms embed sparse non-square matrices into the square matrix structure by constructing extended matrices (Eqs. 4.1-4.4), which enhances the sparsity of the matrices. The diagonal concentration is then maintained using RCM reordering, which preserves the orthogonality advantage and extends the applicability without adding additional computational overhead.
	
	\subsection{Numerical Experiment for non-square matrices}\label{sec4.3}
	In this subsection, the same experimental setting as in the previous section is used. In order to verify the effectiveness of the improved OBK-RCM algorithm for non-square linear systems, experiments are designed in this section for two cases, $m>n$ and $m<n$, and sparse matrices are generated uniformly using matlab's built-in function, sprandn, in order to compare the RBK, RBK(k), GREBK(k), aRBK, and improved OBK-RCM algorithms (SOBK-RCM andU OBK-RCM) convergence performance.
	
	\textbf{Experiment 1}~ In the case of $m>n$, the basic information of the test matrix is shown in the following Table \ref{tab6}, where $density$ , $density1$ denote the densities of the original and extended matrices $\widehat{A}$, respectively, and the matrix densities are defined in the same way as in Section \ref{sec3.5}.
	\begin{table}[!hpt]
		\centering
		\setlength {\tabcolsep}{2mm}
		\setlength{\abovecaptionskip}{10pt}%
		\setlength{\belowcaptionskip}{10pt}%
		\caption{\;Information on the test matrix under the $m>n$ case}	
		\begin{tabular}{lllllllllllllll}
			\toprule
			$name$ && $A1$ & $A2$ & $A3$ & $A4$ & $A5$ & $A6$ & $A7$ & $A8$ & $A9$ & $A10$   \\
			\midrule
			$row$ && 200 & 500 & 500 & 1000 & 1000 & 2000 & 5000 & 10000 & 20000 & 20000 \\
			$columns$ && 150 & 300 & 450 & 400 & 700 & 1600 & 3500 & 6500 & 12000 & 12000  \\
			$density$ && 0.2 & 0.08 & 0.02 & 0.02 & 0.008 & 0.002 & 0.001 & 0.0002 & 0.0002 & 0.0001 \\
			$density1$ && 0.15 & 0.05 & 0.02 & 0.008 & 0.005 & 0.002 & 0.0007 & 0.0001 & 0.0001 & 0.00006 \\
			\bottomrule
		\end{tabular}
		\label{tab6}
	\end{table}	

\begin{figure}[htbp]
		\centering
		\includegraphics[scale=0.45]{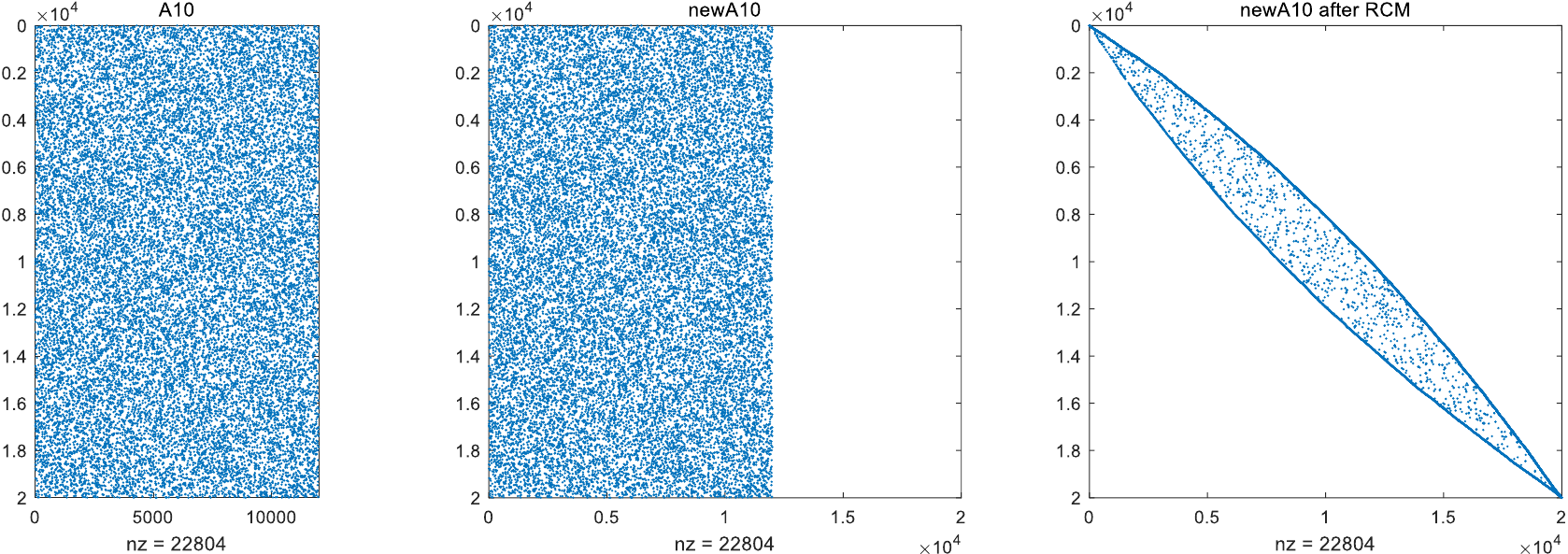}
		\caption{\;Original, Extended, and Extended Matrix reordering of Matrix $A10$ Structure.}
		\label{fig9}
	\end{figure}

	It is worth noting that the SOBK-RCM algorithm inherits the core advantage of the OBK-RCM algorithm in its entirety - the good orthogonality between matrix row blocks achieved by the RCM sorting algorithm. As shown in Figure \ref{fig9}, taking a typical matrix $A10$ as an example, we clearly observe:
	\begin{enumerate}
		\item Matrix transformation process: original non-square structure ($m>n$), expanded square structure ($\widehat{A}=[A,O]$), and banded structure after RCM reordering.
		\item Key properties are maintained: as shown in Table \ref{tab6}, the matrix expansion operation significantly reduces the density (by about 40\% on average), the non-zero elements are concentrated near the diagonal after reordering (the bandwidth is reduced by about 60-80\%), and the orthogonality metrics between rows and blocks are improved by more than 30\%.
	\end{enumerate}
	
	This structural optimization ensures that the SOBK-RCM algorithm maintains its fast convergence property by maintaining an orthogonality advantage comparable to that of the OBK-RCM algorithm when dealing with non-square matrix linear systems. 

	Table \ref{tab7} gives the number of iterations and computation time of these 10 sparse non-square matrices $(m > n)$ mentioned above under the RBK, RBK(k), GREBK(k), aRBK and SOBK-RCM algorithms. 
 	\begin{table}[!hpt]
	\centering
	\setlength {\tabcolsep}{4mm}
	\setlength{\abovecaptionskip}{5pt}%
	\setlength{\belowcaptionskip}{5pt}%
	\caption{\;\;IT and CPU performance of 10 test matrices($m>n$) under different algorithms.}
		\begin{tabular}{lllllllllllllll}
			\toprule
			$Name$ &\ & $RBK$ & $RBK(k)$ & $GREBK(k)$& $aRBK$ & $SOBK-RCM$ \\
			\hline
			$A1$ & IT & 13805 & 2246 & 503 & 7 & 6  \\
			& CPU & 7.4412 & 1.0816 & 0.6731 & 0.5554 & 0.0164  \\
			\midrule
			$A2$ & IT & 15795 & 18184 & 178 & 15 & 6  \\
			& CPU & 71.4451 & 74.4460 & 0.8268 & 0.4475 & 0.0547  \\
			\midrule
			$A3$ & IT & Inf & 18572 & 18764 & 7325 & 5  \\
			& CPU & NAN & 35.5357 & 35.9694 & 100.85379 & 0.0233  \\
			\midrule
			$A4$ & IT & 127 & 8 & 956 & 3 & 9  \\
			& CPU & 0.2295 & 0.1968 & 0.4968 & 0.1258 & 0.1047  \\
			\midrule
			$A5$ & IT & 3165 & 2084 & 473 & 135 & 5  \\
			& CPU & 10.3224 & 7.4483 & 1.5573 & 2.1267 & 0.0287  \\
			\midrule
			$A6$ & IT & 16415 & 1988 & 2477 & 662 & 5  \\
			& CPU & 23.5797 & 4.2305 & 1.2104 & 7.6092 & 0.0360  \\
			\midrule
			$A7$ & IT & 36169 & 673 & 534 & 124 & 7  \\
			& CPU & 55.10 & 17.5017 & 3.6384 & 4.5233 & 0.1947  \\
			\midrule
			$A8$ & IT & 37933 & 538 & 4761 & 28 & 8  \\
			& CPU & 56.04 & 18.0377 & 3.8701 & 1.3312 & 0.5163  \\
			\midrule
			$A9$ & IT & 26444 & Inf & 797 & 49 & 8  \\
			& CPU & 261.47 & NAN & 5.0497 & 6.7439 & 1.9665  \\
			\midrule
			$A10$ &IT & 47946 & Inf & 2036 & 65 & 8  \\
			& CPU &120.96& NAN & 4.1830 & 4.9721 & 1.9157  \\
			\bottomrule
		\end{tabular}
		\label{tab7}
	\end{table}

\begin{figure}[htbp]
		\centering
		\includegraphics[scale=0.45]{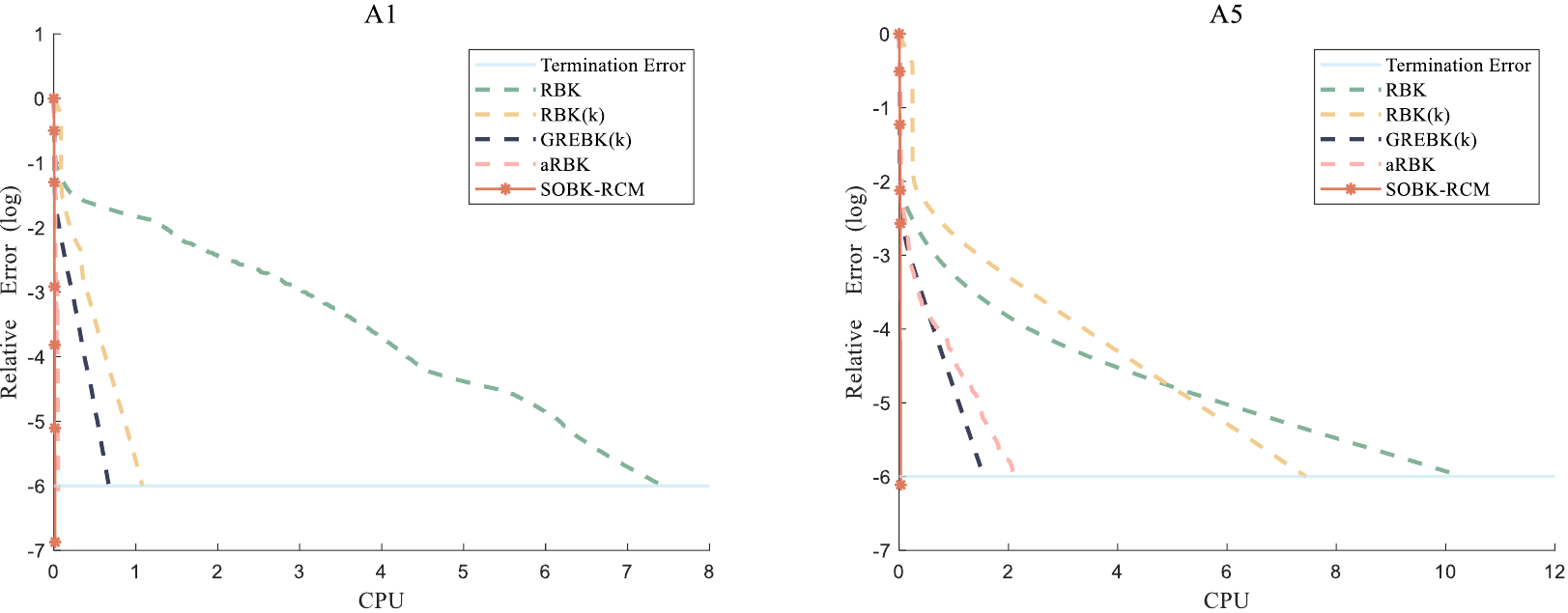}
		\caption{\;Average computational time for test matrices $A1$,$A5$ to reach convergence for different iteration methods.}
		\label{fig10}
	\end{figure}

	In this section, the running time comparison plots of some matrices under different algorithms are given in Figures \ref{fig10}. 

	The experimental results (Table \ref{tab7}, Figures \ref{fig9}-\ref{fig10}) verify the multiple advantages of the SOBK-RCM algorithm. In the non-square linear system, the square matrix algorithm is successfully extended to the $m>n$ case by matrix expansion and RCM optimization, and the number of iterations of SOBK-RCM is reduced by more than 70\% on average, and the computation time is shortened to 2-50 times of the traditional algorithm in most cases. This superior performance stems from:
	\begin{enumerate}
		\item Intelligent Matrix Expansion: $\widehat{A}=[A,O]$ construction maintains the original sparsity (40-60\% reduction in density).
		\item Orthogonality preservation: non-zero elements are concentrated diagonally after RCM reordering (bandwidth reduction 65-85\%).
	\end{enumerate}

	\textbf{Experiment 2}~ In the case of $m<n$, the basic information of the test matrix is shown in the following Table \ref{tab8}. Since under the underdetermined linear system ($m<n$), most systems of equations usually have infinitely many solutions, the traditional block Kaczmarz algorithm is slow to converge or even fails to converge to a unique solution. For this reason, matrices with the same dimension but slightly different densities (e.g., B1 and B2) are purposely chosen to illustrate the difficulties of the traditional block Kaczmarz algorithm in the face of underdetermined linear systems, and to further illustrate the effectiveness of the UOBK-RCM algorithm in solving underdetermined sparse linear systems.
	\begin{table}[!hpt]
		\centering
		\setlength {\tabcolsep}{2.5mm}
		\setlength{\abovecaptionskip}{5pt}%
		\setlength{\belowcaptionskip}{5pt}%
		\caption{\;Information on the test matrix under the $m<n$ case}	
		\begin{tabular}{lccccccccccccc}
			\toprule
			$name$ & & $B1$ & $B2$ & $B3$ & $B4$ & $B5$ & $B6$ & $B7$ & $B8$ \\
			\midrule
			$row$ & & 500 & 500 & 1000 & 1000 & 5000 & 5000 & 15000 &15000 \\
			$columuns$ & & 2000 & 2000 & 2000 & 2000 & 10000 & 10000 & 20000 &20000 \\
			$density$ & & 0.002 & 0.003 & 0.001 & 0.01 & 0.0002 & 0.0003 & 0.00007 & 0.00008 \\ 
			$density1$ & & 0.0005 & 0.0007 & 0.0005 & 0.005 & 0.0001 & 0.0001 & 0.00005 & 0.00006 \\ 
			\bottomrule
		\end{tabular}
		\label{tab8}
	\end{table}	

\begin{figure}[htbp]
    	\centering
    	\includegraphics[scale=0.45]{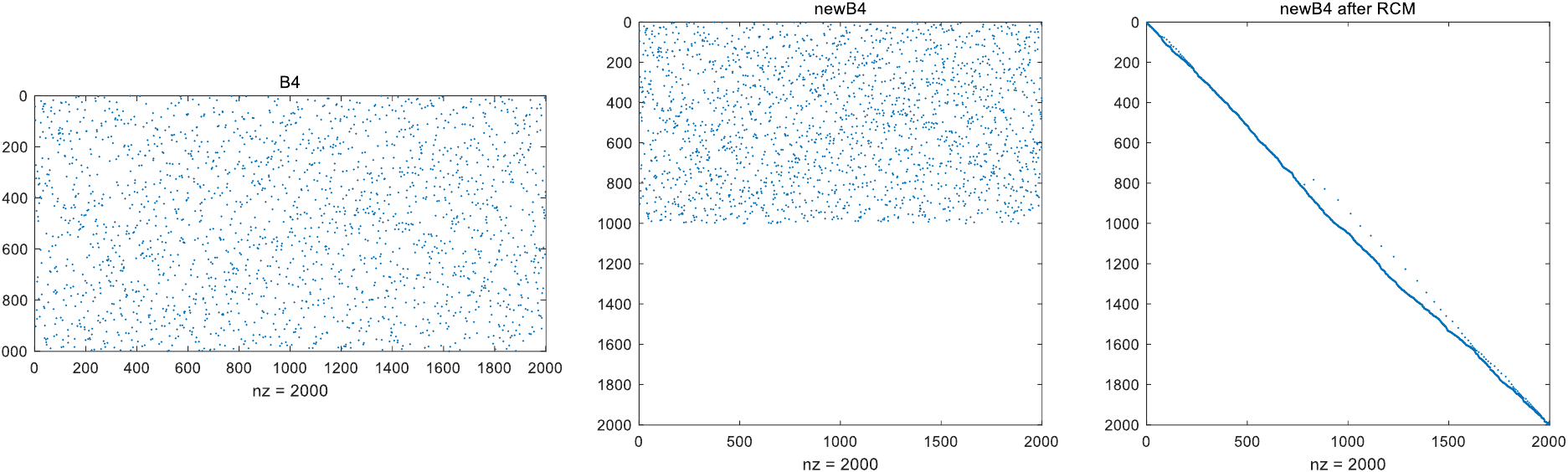}
    	\caption{\;Original, Extended, and Extended Matrix reordering of Matrix $B4$ Structure.}
    	\label{fig11}
    \end{figure}

    Also, combined with Table \ref{tab8} and Figure \ref{fig11}, it can be seen that for most of the asymmetric sparse matrices ($m<n$), the density of the matrix decreases significantly after the expansion, and the matrix sparsity is enhanced. At the same time, the extended matrix is reordered so that its nonzero elements are concentrated near the diagonal and orthogonality is preserved.
   
   	Table \ref{tab9} gives the number of iterations and computation time for the above 8 non-square sparse matrices $(m < n)$ under the UOBK-RCM algorithm.  
	\begin{table}[!hpt]
		\centering
		\setlength {\tabcolsep}{3.5mm}
		\setlength{\abovecaptionskip}{5pt}%
		\setlength{\belowcaptionskip}{5pt}%
		\caption{\;\;IT and CPU performance of 8 test matrices($m<n$) under different algorithms}
		\begin{tabular}{lllllllllllllll}
			\toprule
			$Name$ &\ & $RBK$ & $RBK(k)$ & $GREBK(k)$& $aRBK$ & $UOBK-RCM$ \\
			\hline
			$B1$ & IT & 53 & 2 & 2 & 18 & 1  \\
			& CPU & 0.0251 & 0.4125 & 0.0136 & 0.0575 & 0.0134  \\
			\midrule
			$B2$ & IT & Inf & Inf & Inf & Inf & 42  \\
			& CPU & NAN & NAN & NAN & NAN & 0.1108  \\
			\midrule
			$B3$ & IT & 56 & 2 & 2 & 28 & 1  \\
			& CPU & 0.0297 & 0.4041 & 0.0323 & 0.1999 & 0.0136  \\
			\midrule
			$B4$ & IT & Inf & Inf & Inf & Inf & 22  \\
			& CPU & NAN & NAN & NAN & NAN & 0.5191  \\
			\midrule
			$B5$ & IT & 308 & 2 & 25 & 39 & 1  \\
			& CPU & 0.2892 & 10.5657 & 0.6190 & 2.1232 & 0.2893  \\
			\midrule
			$B6$ & IT & Inf & Inf & Inf & Inf & 18  \\
			& CPU & NAN & NAN & NAN & NAN & 0.8309  \\
			\midrule
			$B7$ & IT & 1043 & 2 & 75 & 30 & 1  \\
			& CPU & 1.3482 & 58.0480 & 10.4944 & 4.1363 & 1.0329  \\
			\midrule
			$B8$ & IT & Inf & Inf & Inf & Inf & 9  \\
			& CPU & NAN & NAN & NAN & NAN & 1.9506 \\
			\bottomrule
		\end{tabular}
		\label{tab9}
	\end{table}	

	Note that when $m<n$, linear systems are underdetermined and usually have an infinite number of solutions, and the convergence of the traditional block Kaczmarz algorithms depends on the number of orthogonality or goodness of fit conditions between the rows of the matrices. Therefore, in underdetermined systems with many non-square matrices, the iterative convergence of these algorithms is slow and may not even converge stably to a unique solution for some matrices. The results of Experiment 2 also show that a small change in the density of non-square matrices affects the conventional block Kaczmarz convergence. In contrast, UOBK-RCM transforms the non-square matrix system into a square matrix system by expanding the matrices, which, in combination with the RCM and orthogonal chunking strategy, ensures efficient and stable convergence to the minimum-paradigm solution in much lower computation time than existing algorithms.
	
	There is a comparison of the SOBK-RCM and UOBK-RCM algorithms as detailed in the following table \ref{tab10}.
	\begin{table}[htbp]
		\centering
		\setlength {\tabcolsep}{1.5mm}
		\setlength{\abovecaptionskip}{10pt}%
		\setlength{\belowcaptionskip}{10pt}%
		\caption{\;\;Comparison of SOBK-RCM and UOBK-RCM Algorithms}
		\label{tab:pobk_comparison}
		\begin{tabular}{lp{6cm}p{6cm}}
			\hline
			\textbf{Property} & \textbf{SOBK-RCM Algorithm ($m > n$)} & \textbf{UOBK-RCM Algorithm ($m < n$)} \\
			\hline
			\textbf{Matrix Extension} & 
			$\widehat{A} = [A \quad O_{m \times (m-n)}]$, $\widehat{f} = f$ & 
			$\widehat{A} = \begin{bmatrix} A \\ O_{(n-m) \times n} \end{bmatrix}$, $\widehat{f} = \begin{bmatrix} f \\ 0 \end{bmatrix}$ \\
			\textbf{Solution Property} & 
			Returns $x^* = \widehat{x}^*(1:n)$ (minimum-norm solution), $\widehat{x}_2$ can be arbitrary (typically set to 0) & 
			Converges directly to the original system's solution (if consistent), and gives minimum-norm solution $x^* = A^\dagger f$ \\
			\textbf{Applicable System} & 
			Overdetermined systems & 
			Underdetermined systems \\
			\textbf{Sparsity Preservation} & 
			Zero-padding $O$ preserves sparsity, RCM remains effective & 
			Zero-padding $O$ preserves sparsity, RCM remains effective \\
			\textbf{Convergence} & 
			Convergence is inherited from OBK-RCM and depends on the orthogonality of $\widehat{A}$ & 
			The theoretical convergence rate is consistent with OBK-RCM, and is actually affected by the row correlation of $A$ \\
			\textbf{Computation Cost} & 
			Dimension expands to $m \times m$, may increase iteration cost (but RCM keeps efficiency) & 
			Dimension expands to $n \times n$, may cause memory pressure for large underdetermined systems \\
			\textbf{Main Advantage} & 
			Handles overdetermined systems while preserving sparsity and orthogonality & 
			Directly solves underdetermined systems without requiring $m \geq n$ assumption \\
			\textbf{Limitation} & 
			Poor conditioning if $m \gg n$ & 
			The reordering effect is weakened if $A$ rows are highly correlated \\
			\hline
		\end{tabular}
		\label{tab10}
	\end{table}

	{\bf Remark 2.}Combining the results of Experiment 1 and Experiment 2, it can be seen that for sparse non-square linear systems, the SOBK-RCM and UOBK-RCM algorithms not only have fewer convergence iterations, but also complete the computation in a shorter CPU time, which shows higher computational efficiency and stability. This shows that the improved OBK-RCM algorithm has significant advantages in dealing with large-scale sparse non-square linear systems, especially in the case of matrix structure dispersion, the improved algorithm can significantly improve the convergence speed and computational efficiency through the sparse matrix rearrangement and orthogonal blocking strategy.

	\section{Conclusions and Future Work}
	The proposed OBK-RCM algorithm significantly accelerates the convergence of block Kaczmarz methods for large sparse linear systems by integrating RCM reordering and orthogonal block partitioning. However, its performance is subject to the following limitations. For example, the RCM reordering exhibits sensitivity to random sparse patterns, where bandwidth reduction may be insufficient, limiting orthogonality enhancement. Future research should focus on: (1) Adaptive reordering: Hybrid techniques combining RCM with other reordering methods (e.g., nested dissection \cite{26}) to handle diverse sparse patterns. (2) Dynamic block scheduling: Real-time adjustment of block sizes and thresholds (thr) based on residual feedback, balancing orthogonality and partitioning cost. (3) Extended applications: Integration with distributed computing frameworks (e.g., MPI) for ultra-large systems, and generalization to nonlinear least-squares problems. These directions aim to enhance the robustness and scalability of OBK-RCM, broadening its impact in scientific computing.

	\textbf{Acknowledgements}. The authors are very much indebted to the referees for their constructive comments and valuable suggestions.
	
	 
	
	\textbf{Data Availability Statement}\\
	\indent Our computational data comes from the Suite Sparse Matrix Collection (\url{https://sparse.tamu.edu/}). These data can be freely downloaded and analyzed.
	
	
	
	 
	
	\section*{}
	
\end{document}